\documentclass[12pt]{elsart}
\usepackage{amsfonts}
\usepackage{epsf}
\journal{Physica D}

\newcommand{\be}{\begin{equation}}
\newcommand{\ee}{\end{equation}}
\newcommand{\bea}{\begin{eqnarray}}
\newcommand{\eea}{\end{eqnarray}}
\newcommand{\beaa}{\begin{eqnarray*}}
\newcommand{\eeaa}{\end{eqnarray*}}
\newcommand{\eq}[1]{(\ref{#1})}

\newcommand{\Rl}{\mathbb{R}}
\newcommand{\Ir}{\mathbb{Z}}
\newcommand{\Sp}{\mathbb{S}}
\newcommand{\RP}{\mathbb{R}P}
\newcommand{\eps}{\varepsilon}
\renewcommand{\phi}{\varphi}

\renewcommand{\div}{\mathrm{div}\,}
\newcommand{\curl}{\mathrm{curl}\,}

\def\bn{\mathbf{n}}
\def\bB{\mathbf{B}}
\def\bk{\mathbf{k}}
\def\bx{\mathbf{x}}
\def\a{A}

\def\bm{\mathbf{m}}
\def\bu{\mathbf{u}}
\def\bF{\mathbf{F}}

\def\bnz{\mathbf{n}_0}

\def\bnp{\mathbf{n}^\prime}
\def\ap{A^\prime}
\def\bBp{\mathbf{B}^\prime}
\def\lambdap{\lambda^\prime}

\def\bnh{\widehat{\mathbf{n}}}

\def\lambdah{\widehat{\lambda}}

\def\tA{\widehat{A}}
\def\tB{\widehat{\mathbf{B}}}

\def\rA{\widetilde{A}}
\def\rB{\widetilde{\mathbf{B}}}
\def\rlambda{\widetilde{\lambda}}

\def\ba{\mathbf{a}}
\def\bb{\mathbf{b}}

\def\Wtwist{b}
\def\Wsplay{a}
\def\Rtwist{\mathbf{S}}
\def\Rsplay{\mathbf{R}}

\newtheorem{proposition}{Proposition}

\begin{document}
\begin{frontmatter}

\title{Orientation Waves in a Director Field\\ with Rotational Inertia}

\author[ga]{Giuseppe Al\`\i},
\ead{g.ali@iac.cnr.it}
\author[jkh]{John K. Hunter\corauthref{cor}}
\corauth[cor]{Corresponding author}
\ead{jkhunter@ucdavis.edu}

\address[ga]{Istituto per le Applicazioni del Calcolo \\
Consiglio Nazionale delle Ricerche, Napoli \\
and INFN-Gruppo c. Cosenza, Italy}

\address[jkh]{Department of Mathematics \\
University of California at Davis\\
Davis CA 95616, USA}

\date{September 6, 2006}

\maketitle

\begin{abstract}
We study the propagation of orientation waves in a director field
with rotational inertia and potential energy given by the Oseen-Frank
energy functional from the continuum theory of nematic liquid crystals.
There are two types of waves,
which we call splay and twist waves.
Weakly nonlinear splay waves are described
by the quadratically nonlinear Hunter-Saxton equation.
Here, we show that weakly nonlinear twist waves
are described by a new cubically nonlinear, completely
integrable asymptotic equation. This equation provides
a surprising representation of the Hunter-Saxton equation
for $u$ as an advection equation for $v$, where
$u_{xx} = v_x^2$. There is
an analogous representation of the Camassa-Holm equation.
We use the asymptotic equation to analyze
a one-dimensional initial value problem for the
director-field equations with twist-wave initial data.
\end{abstract}

\end{frontmatter}

\section{Introduction}
\setcounter{equation}{0}

In this paper, we analyze a system of partial differential equations
that models
the propagation of orientation waves in a massive director field.
The restoring force for the waves
is provided by the Oseen-Frank energy \eq{o-f}, used in the
continuum theory of nematic liquid crystals,
and the inertia is provided by the rotational inertia
of the director field.

The motion of the director field is then described by the variational
principle \eq{varprinc}. The corresponding Euler-Lagrange equation \eq{dirPDE} is scale-invariant,
and forms a non-dispersive hyperbolic system of wave equations. The system
supports two types of waves, which we call `splay' and `twist' waves, respectively.
The principal nonlinear effects are that the wave speeds
depend on the angle between the propagation
direction and the director field, and that the twist waves generate splay waves.

Splay waves were investigated by Saxton \cite{saxton} and Hunter and Saxton \cite{hs}. In that case, the
one-dimensional system of equations for the director field (see Section~\ref{sec:oned})
reduces to a scalar wave equation for an angle $\phi(x,t)$,
\be
\phi_{tt} - \left[a^2(\phi)\phi_x\right]_x + a(\phi) a^\prime(\phi) \phi_x^2 = 0.
\label{scalarwaveeq}
\ee
Here, the wave-speed
$a$ is a smooth, nonzero function of $\phi$, given in \eq{dircfs},
and the prime denotes
a derivative with respect to $\phi$.

The wave equation \eq{scalarwaveeq} is obtained from the variational principle
\[
\delta \int \frac{1}{2}\left\{\phi_t^2 - a^2(\phi) \phi_x^2\right\}\,dx dt = 0,
\]
which is one of the simplest nonlinear
generalizations of the variational principle for the linear wave equation
one can imagine. The effects of this nonlinearity
include the formation of cusp-type singularities \cite{ghz}.
Smooth solutions can be extended by global weak solutions,
but in sharp contrast with the more familiar case of entropy solutions
of hyperbolic conservation laws \cite{dafermos}, equation \eq{scalarwaveeq}
possesses conservative weak solutions that are compatible
with its variational and Hamiltonian structure. The global well-posedness
of the initial value problem for \eq{scalarwaveeq}
for conservative weak solutions is proved in \cite{bz}.

The use of weakly nonlinear asymptotics to study the
behavior of solutions of \eq{scalarwaveeq} that consist of small, localized
perturbations
$u(x,t)$ of a constant state $\phi_0$ leads (after an appropriate normalization,
and provided that $a^\prime(\phi_0)\ne 0$)
to the Hunter-Saxton (HS) equation \cite{hs},
\be
\left(u_t + u u_x\right)_x - \frac{1}{2} u_x^2 = 0.
\label{HSeq}
\ee
This equation is derived from
the variational principle
\be
\delta \int \frac{1}{2}\left\{u_x u_t + u u_x^2 \right\} \, dx dt = 0.
\label{hsaction}
\ee
It is completely integrable \cite{hz,bs, km}, and possesses global
dissipative and conservative weak solutions \cite{hz1,zz,bc}.

In this paper, we study the full system of director-field equations
and show that qualitatively new phenomena arise
which are not present for the scalar wave equation \eq{scalarwaveeq}.
The structure of the system is most easily seen in the case
of one space-dimension; the director-field equations then reduce
to two coupled wave equations \eq{wavea}--\eq{waveb}
for angles $\phi(x,t)$, $\psi(x,t)$,
given by the variational principle \eq{onedvar}.
The angle $\phi$, with wave-speed $a$,
corresponds to splay waves; the angle $\psi$, with wave speed
$b$, corresponds to twist waves. Both wave speeds $a$, $b$ depend only
on $\phi$, and the wave equation for $\phi$ is forced by
source terms that are proportional to quadratic
functions of the derivatives of $\psi$.
If $\psi$ is constant, then the system \eq{wavea}--\eq{waveb} reduces to
the scalar wave equation \eq{scalarwaveeq}.

Our main result is the following asymptotic PDE for a localized,
weakly-nonlinear twist wave with amplitude $v(x,t)$:
\bea
&&\left(v_t + u v_x\right)_x = 0,
\label{uveq}
\\
&&u_{xx} = v_x^2.
\nonumber
\eea
This cubically-nonlinear, scale-invariant system
consists of an advection equation for $v$ in which the
advection velocity $u$ is reconstructed
nonlocally and quadratically from $v$.
It is obtained from the variational principle
\be
\delta \int \frac{1}{2}\left\{v_x v_t + u v_x^2 + \frac{1}{2} u_x^2\right\} \, dx dt = 0.
\label{action}
\ee
The second dependent variable $u(x,t)$ in \eq{uveq} may be interpreted as
the amplitude of a splay wave that is generated nonlinearly by the twist wave, which
then affects the velocity of the twist wave.

Remarkably, the elimination of $v$ from \eq{uveq} implies that $u$ satisfies
the HS-equation
\be
\left[ \left(u_t + u u_x\right)_x - \frac{1}{2} u_x^2\right]_x = 0.
\label{HSeq1}
\ee
Thus, a solution $v$ of \eq{uveq}
is related to a solution $u$ of \eq{HSeq1} by $v_x^2 = u_{xx}$.
Under this change of variables,
the variational principle \eq{action} transforms into
a second variational principle
for \eq{HSeq1}, distinct from \eq{hsaction} but compatible with it,
that is associated with a Lie-Poisson structure (see Proposition~\ref{prop:hamtrans}).

The HS-equation therefore arises from the system of director-field
equations in two different asymptotic limits for two different types of waves,
with two separate variational structures.
The resulting
bi-Hamiltonian structure explains why \eq{HSeq1}
is completely integrable \cite{hz}, and it
follows that \eq{uveq} is also completely integrable (see Section~\ref{sec:ham}).
Furthermore, as we show in Section~\ref{sec:char},
equation \eq{uveq} may be solved explicitly by the method of characteristics; the key point
is that the Jacobian of the transformation between spatial and characteristic coordinates
satisfies an integrable Liouville equation.

The correspondence between solutions of  \eq{uveq} and \eq{HSeq1}
is neither one-to-one or onto.
Only convex solutions of \eq{HSeq1} with $u_{xx} \ge 0$
can be obtained from solutions of \eq{uveq}, and
solutions $v$ of \eq{uveq} for which $v_x$
has the same magnitude but different signs (which may depend on $x$)
correspond to the same solution of \eq{HSeq1}.
Moreover, because of the nonlinear nature of the transformation,
distributional solutions of \eq{HSeq1}
do not necessarily transform into distributional solutions of \eq{uveq}.
Nevertheless, equation \eq{uveq} provides
an interesting representation of the HS-equation \eq{HSeq1} for $u$ as an advection equation
for a new variable $v$. There is a similar representation of the
closely related Camassa-Holm (CH) equation \cite{CH} (see \eq{sys1}--\eq{sys2} below,
with $Mu = u_{xx}-u$).


The twist and splay waves exhibit a basic difference in their nonlinearity.
Twist waves are linearly degenerate, meaning that the derivative of their
wave speed with respect to the wave amplitude is identically zero, as stated in \eq{lindeg}.
This linear degeneracy may be seen
in the asymptotic system \eq{uveq}, where the advection velocity $u$ of $v$ is
independent of $v$, and in the one-dimensional
system \eq{wavea}--\eq{waveb}, where the wave speed $b(\phi)$ of $\psi$
is independent of $\psi$.
By contrast, excluding exceptional values of $\phi$ where $a^\prime(\phi) = 0$,
splay waves are genuinely nonlinear, meaning that the derivative of their
wave speed with respect to the wave amplitude is non-zero, as stated in \eq{gnl}.
Despite the linear degeneracy of the twist waves, their interaction with
splay waves leads to nontrivial, cubically nonlinear dynamics.

Equation \eq{uveq} differs from another, more obvious, cubically-nonlinear,
scale-invariant modification of the HS-equation \cite{hs,zz1},
\be
\left(u_t + u^2 u_x\right)_x = u u_x^2,
\label{cubic1}
\ee
given by the variational principle
\[
\delta \int \frac{1}{2}\left\{u_x u_t + u^2 u_x^2 \right\} \, dx dt = 0.
\]
Unlike \eq{uveq}, equation \eq{cubic1} is an asymptotic limit of the scalar wave equation \eq{scalarwaveeq},
and arises when there is a loss of genuine nonlinearity at the unperturbed state $\phi_0$ such that
$a^\prime(\phi_0) = 0$ but $a^{\prime\prime}(\phi_0) \ne 0$. In Section~\ref{sec:rot},
we derive from the full director-field system a vector generalization of \eq{cubic1},
given in \eq{normpoleq} or \eq{pol1}--\eq{pol2},
that describes the propagation of polarized orientation waves
in the same direction as the unperturbed director field.

We now outline the contents of this paper.
In Section~\ref{sec:direq},
we describe the system of PDEs for the director-field,
and briefly compare it with some related variational field theories.
We show that the system is hyperbolic, and study
the linearized orientation waves. In Section~\ref{sec:wn},
we summarize the asymptotic equations for
weakly nonlinear splay and twist waves.

In Section~\ref{sec:oned}, we
write out the one-dimensional director-field equations in terms of
spherical polar angles, and
construct an asymptotic solution of a one-dimensional initial value problem
with initial data corresponding to a small-amplitude, compactly supported twist wave.
This solution illustrates, in particular, the generation of splay waves by twist waves,
and our main goal is to formulate equations for the resulting splay waves.

There are two cases, depending on whether the twist waves are faster or slower
than the splay waves (see Figure~\ref{charivp}). When the twist waves are faster,
they move at a constant velocity into an unperturbed director field,
and we obtain a Cauchy problem
for the splay-wave equation \eq{scalarwaveeq} with data for $\phi$, $\phi_x$
given on the space-like twist-wave trajectories. When the twist waves are slower,
they are embedded inside the splay wave, and
we obtain a free-boundary problem for \eq{scalarwaveeq}.
The twist-wave trajectories are not known \textit{a priori}, and
$\phi_x$ satisfies a jump condition
across the trajectories that follows from an integration of
the asymptotic equation \eq{uveq}.

In Section~\ref{sec:char}, we solve \eq{uveq} by the method of characteristics,
and use the result to solve the initial-boundary value problem
for \eq{uveq} that arises in Section~\ref{sec:oned}.
In Section~\ref{sec:ham}, we show that \eq{uveq} leads to
the HS-equation, and describe its bi-Hamiltonian structure
and integrability. We also consider the CH-equation.
In Section~\ref{sec:per}, we derive an apparently non-integrable
generalization of \eq{uveq}, given in \eq{uvpereq1}--\eq{uvpereq2},
that applies to periodic twist waves with mean-field interactions.


\section{Director fields}
\label{sec:direq}
\setcounter{equation}{0}

The system of wave equations we study here is motivated by the theory of nematic
liquid crystals. We consider an unbounded anisotropic medium whose orientation (in three space dimensions)
is described by a time-dependent director field $\bn(\bx,t)$ of unit vectors,
$\bn : \Rl^3\times \Rl \to \Sp^2$.

For uniaxial nematic liquid crystals, the directions $\bn$ and $-\bn$
are identified, and then $\bn : \Rl^3\times \Rl \to \RP^2$. In most of the problems we consider
here, the director field is a small perturbation of a given smooth field,
and the global topological distinctions between $\Sp^2$ and $\RP^2$ are not relevant.

We suppose that the potential energy $W$ of a director field $\bn$ is given by
the Oseen-Frank energy functional \cite{deG}
\be
W\left(\bn,\nabla\bn\right) = \frac{1}{2}\alpha \left(\div\bn\right)^2
+ \frac{1}{2} \beta\left(\bn\cdot\curl \bn\right)^2
+ \frac{1}{2} \gamma\left|\bn\times\curl \bn\right|^2.
\label{o-f}
\ee
Up to a null-Lagrangian, $W$ is the most general quadratic function of
$\nabla \bn$ with coefficients depending on $\bn$ that is invariant under the transformations
$\bn\mapsto -\bn$ and
\be
\bx \mapsto R\bx,\qquad \bn \mapsto R\bn\qquad\mbox{for all orthogonal maps $R$}.
\label{lqinv}
\ee
The positive coefficients $\alpha$, $\beta$, $\gamma$ are elastic constants
of splay, twist, and bend, respectively. For liquid crystals, we typically
have $0 < \beta < \alpha < \gamma$.

In the one-constant approximation,
$\alpha = \beta = \gamma$,
the Oseen-Frank energy reduces (up to
a null-Lagrangian) to the harmonic map energy,
\[
W(\nabla\bn) = \frac{1}{2}\alpha \left|\nabla\bn\right|^2.
\]
This energy does not depend explicitly on $\bn$, and is
invariant under a larger group of transformations,
namely
\be
\bx \mapsto R\bx,\qquad \bn \mapsto S\bn\qquad\mbox{for all orthogonal maps $R$, $S$}.
\label{hminv}
\ee
The nonlinear effects we study here
vanish in this case, and we assume that
the elastic constant are distinct.

We further suppose that the director field has rotational kinetic energy
and is not subject to damping or dissipation.
Its motion is then governed by the
variational principle
\be
\delta \int_{\Rl^3\times\Rl} \left\{\frac{1}{2} \bn_t^2 - W\left(\bn,\nabla\bn\right)\right\}
\, d\bx dt = 0,
\qquad \bn\cdot\bn = 1,
\label{varprinc}
\ee
where we have normalized the moment of inertia per unit volume
of the director field to one.

The Euler-Lagrange equation associated with \eq{varprinc} is
\bea
&&\bn_{tt} = \alpha\nabla(\div \bn) - \beta\left\{\a\curl\bn + \curl(\a\bn)\right\}
\nonumber\\
&&\qquad\qquad\quad
+ \gamma\left\{\bB\times\curl\bn - \curl\left(\bB\times\bn\right)\right\} + \lambda\bn,
\label{dirPDE}
\eea
where
\[
\a = \bn\cdot\curl\bn,\qquad
\bB = \bn\times \curl\bn.
\]
The Lagrange multiplier $\lambda(\bx,t)$ is chosen so that
$\bn\cdot \bn = 1$,
which implies that
\[
\lambda = -\bn_t^2 + \alpha\left[|\nabla \bn|^2 - |\curl \bn|^2\right]
+ 2\left[\beta A^2 + \gamma \bB^2\right] + (\alpha-\gamma)\div\bB.
\]

Equation \eq{dirPDE} provides a natural, geometrical model for
the propagation of orientation waves in an anisotropic medium,
and is representative of a large class
of variational systems of wave equations in which the
wave speeds are functions of the dependent variables \cite{ah}.
It is not applicable to standard liquid crystal
hydrodynamics where the motion of the director field is dominated
by viscosity and the effects of rotational inertia are negligible.
Nevertheless, it is conceivable that the nonlinear
phenomena analyzed here could be observed in high-frequency excitations of
liquid crystals whose molecules possess a large
moment of inertia (such as, perhaps,
liquid crystals obtained from carbon nanotubes).

At a more general level, liquid crystalline phases are the result of
a continuous breaking of rotational symmetry, and the orientation waves we analyze here
may be regarded as associated Goldstone modes \cite{cmp}. Similar phenomena should
occur for non-dispersive Goldstone modes in other systems
with continuously broken symmetry in which the energy density
associated with a set of order parameters $\Psi^a$ has the anisotropic form
\[
A_{ab}^{ij}(\Psi^c) \frac{\partial \Psi^a}{\partial x^i}\frac{\partial \Psi^b}{\partial x^j},
\]
rather than the more commonly assumed isotropic Ginzburg-Landau form proportional to
$|\nabla \Psi|^2$. For liquid crystals, the isotropic form arises only in
the one-constant approximation.

As liquid crystals illustrate, the Ginzburg-Landau energy density is not
dictated by general symmetry arguments in anisotropic media.
For long-wave variations, it is natural to retain the leading-order
terms in the energy that are quadratic in the spatial derivatives of the order parameters,
but the order parameters themselves may vary by a large amount, in which case one
should retain any dependence of the coefficients of the spatial derivatives on the
order parameters.

A similar situation occurs in classical field theories, where nonlinear sigma-models
lead to wave equations with wave speeds that are independent
of the dependent variables. This class of field theories includes the wave-map equations
from Minkowski space into $\Sp^2$ to which \eq{dirPDE} reduces in the one-constant approximation.
On the other hand, general relativity leads to a form of nonlinearity
that is analogous to that of the general director-field equations:
the wave operator in the Einstein equations acts on the metric
and has coefficients that are functions of the metric.
The effects of this nonlinearity in the Einstein equations are, however,
much more degenerate than in the director-field equations \cite{ah1}.

We remark that the invariance properties of these
equations also differ in an analogous way. Nonlinear sigma-models are
invariant under separate transformations of the independent and dependent variables,
as in \eq{hminv}, whereas the gauge-invariance of the Einstein equations
involves a simultaneous transformation of the independent and dependent variables,
as in \eq{lqinv}.

\subsection{Orientation waves}

In this section, we show that \eq{dirPDE} forms a hyperbolic system
of PDEs, and describe the corresponding waves.
We consider solutions of \eq{dirPDE} of the form
\[
\bn(\bx,t) = \bnz + \bnp(\bx,t),
\]
where $\bnp$ is a small perturbation of a constant director
field $\bnz$, linearize the resulting equations for $\bnp$,
and look for Fourier solutions of the linearized equations of the form
\[
\bnp(\bx,t) = N e^{i\bk\cdot\bx - i\omega t} \bnh.
\]
Here, $N\in \Rl$ is an arbitrary amplitude, $\bk\in \Rl^3$ is the wavenumber vector,
$\omega \in \Rl$ is the frequency,
and $\bnh\in \Rl^3$ is a normalized constant vector.

\begin{figure}
\epsfxsize=3in
\epsfysize=3in
\centerline{\epsfbox{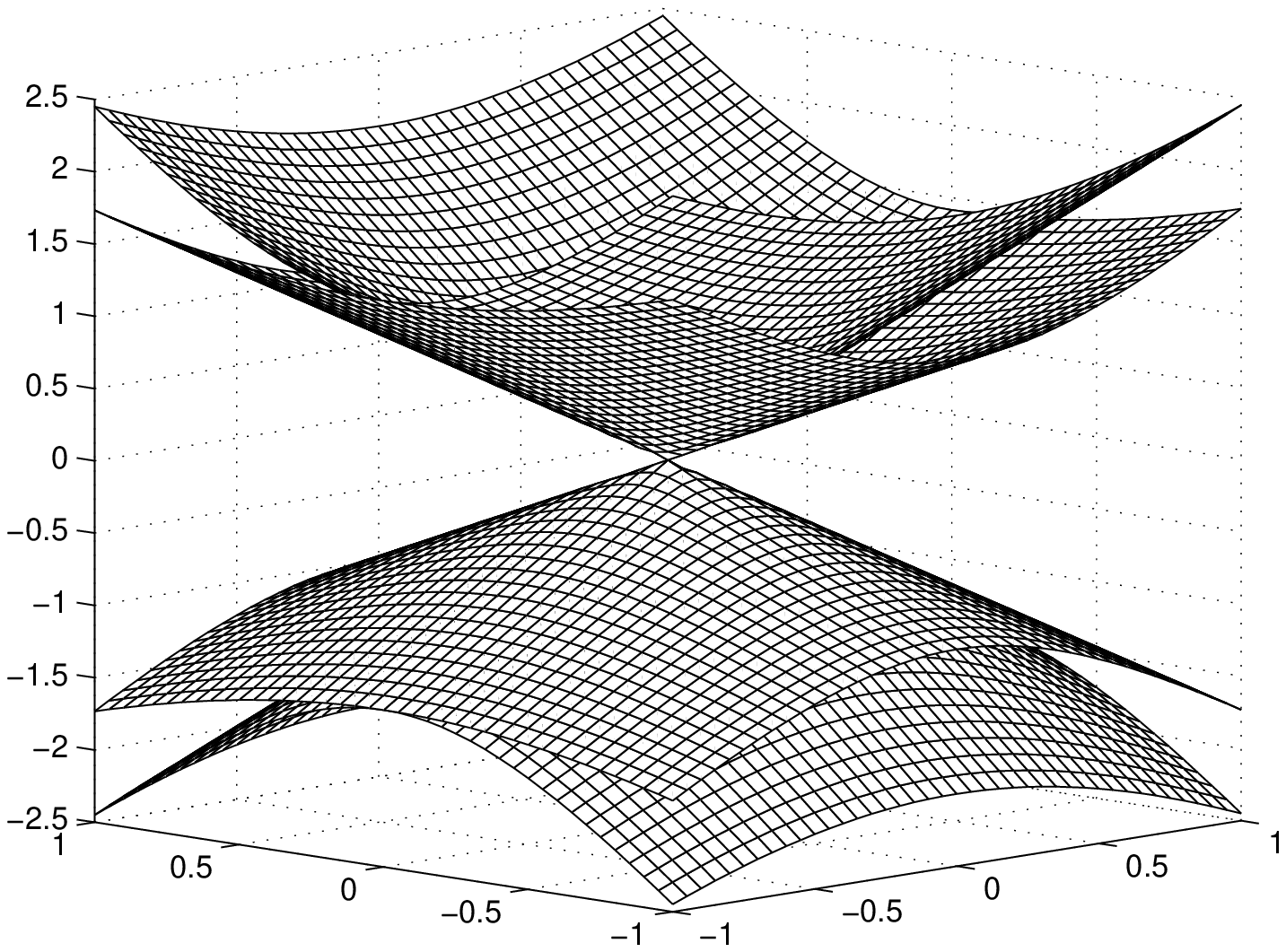}}
\caption{Characteristic variety of \eq{dirPDE}. \label{charvar}}
\end{figure}

We find that (see Appendix \ref{app:lin}) $\omega$ satisfies the linearized dispersion relation
\be
\left[\omega^2 - \Wsplay^2(\bk;\bnz)\right]\left[\omega^2 - \Wtwist^2(\bk;\bnz)\right] = 0,
\label{lindisp}
\ee
where, with $k = |\bk|$,
\bea
&&\Wsplay^2(\bk;\bnz) = \alpha\left[k^2 - \left(\bk\cdot\bnz\right)^2\right] + \gamma\left(\bk\cdot\bnz\right)^2,
\label{amode}
\\
&&\Wtwist^2(\bk;\bnz) = \beta\left[k^2 - \left(\bk\cdot\bnz\right)^2\right] + \gamma\left(\bk\cdot\bnz\right)^2.
\label{bmode}
\eea
Thus, if $\alpha\ne \beta$, the characteristic variety of \eq{dirPDE} consists of two
nested elliptical cones (see Figure~\ref{charvar}). There is a loss of strict hyperbolicity
when $\bk$ is parallel to $\bnz$, corresponding to a wave that propagates in the same direction
as the unperturbed director field.

The waves associated with the branch $\omega^2 = a^2(\bk;\bnz)$
carry perturbations of the director field in the
same plane as $\bnz$ and $\bk$, in which $\bnh = \Rsplay$ where
\be
\Rsplay\left(\bk;\bnz\right) = \bk -  \left(\bk\cdot \bnz\right) \bnz.
\label{Rsplay}
\ee
We call these waves \emph{splay} waves. We note that $\Rsplay$
is orthogonal to $\bnz$, as required by the linearization of the constraint
that $\bn$ is a unit vector.

The waves associated with the branch $\omega^2 = b^2(\bk;\bnz)$
carry transverse perturbations of the director field orthogonal to
$\bnz$ and $\bk$, in which $\bnh = \Rtwist$ where
\be
\Rtwist\left(\bk;\bnz\right) = \bk\times\bnz.
\label{Rtwist}
\ee
We call these waves \emph{twist} waves.

There is a fundamental difference in the dependence of the splay and
twist waves on $\bnz$.
From \eq{amode}, we compute that the derivative of the splay-wave speed
in the direction of the perturbation carried by the wave is given by
\be
\nabla_{\bnz} \Wsplay(\bk;\bnz)\cdot \Rsplay(\bk;\bnz) = -\frac{\left(\alpha-\gamma\right)}{\Wsplay(\bk;\bnz)}
\left(\bk\cdot\bnz\right)\left[k^2 - \left(\bk\cdot\bnz\right)^2\right].
\label{gnl}
\ee
If $\alpha\ne \gamma$, this quantity is nonzero provided that $\bk$ is not
parallel or orthogonal to $\bnz$. On the other hand,
from \eq{bmode}, we see that the derivative of the twist-wave speed
in the direction of the wave is identically zero,
\be
\nabla_{\bnz} \Wtwist(\bk;\bnz)\cdot \Rtwist(\bk;\bnz) = 0.
\label{lindeg}
\ee

By analogy with the terms introduced by Lax in the context of first-order hyperbolic
systems of conservation laws \cite{dafermos}, we say that a wave is \emph{genuinely nonlinear}
if the derivative of the wave speed in the direction
of the wave is nonzero, and \emph{linearly degenerate} if the derivative is identically zero.
Thus, the splay waves are genuinely nonlinear when they do not propagate in directions
parallel
or orthogonal to the director field,
and the twist waves are linearly degenerate.

\section{Weakly nonlinear waves}
\label{sec:wn}
\setcounter{equation}{0}

In this section, we describe the asymptotic equations for
weakly nonlinear splay and twist waves. We also consider the case
of waves that propagate in the same direction
as the unperturbed director field,
when there is a loss of strict hyperbolicity and the
orientation waves are polarized. The algebraic details of the
derivations are summarized in Appendix~\ref{app:details}.

\subsection{Splay waves}
\label{sec:splay}

Weakly nonlinear asymptotics for genuinely nonlinear splay waves leads to the
quadratically nonlinear Hunter-Saxton (HS) equation \cite{hs}. We summarize here the
expansion for weakly nonlinear, non-planar splay waves.

We look for a high-frequency asymptotic solution $\bn^\eps$ of \eq{dirPDE}
with phase $\Phi(\bx,t)$, depending
on a small parameter $\eps$, of the form
\bea
&&\bn^\eps(\bx,t) = \bn\left(\frac{\Phi(\bx,t)}{\eps},\bx,t;\eps\right),
\label{nsplayexp}
\\
&&\bn(\theta,\bx,t;\eps) = \bnz(\bx,t) + \eps \bn_1(\theta,\bx,t) + O(\eps^2)
\qquad\mbox{as $\eps \to 0^+$}.
\label{nsplayexp1}
\eea
In \eq{nsplayexp1}, $\bnz(\bx,t)$ is a smooth solution of \eq{dirPDE} that is independent of $\eps$.

We consider localized waves, such as pulses or fronts,
rather than periodic waves. The resulting
asymptotic solutions are valid near the wavefront
$\Phi(\bx, t)=0$, where $\theta = O(1)$, but they
need not be uniformly valid in $\theta$. Global solutions may be obtained by
matching these `inner' solutions for localized waves with  suitable
`outer' solutions (see Section~\ref{sec:oned}). For periodic waves, there are
additional mean-field interactions (see Section~\ref{sec:per}).

We define the local frequency $\omega(\bx,t)$ and wavenumber $\bk(\bx,t)$ by
\be
\omega = -\Phi_t,
\qquad
\bk = \nabla \Phi.
\label{locfreq}
\ee
We assume that $\bk\ne 0$.
It follows from the expansion that $\omega$, $\bk$
satisfy the linearized dispersion relation \eq{lindisp},
and we suppose that they satisfy the
splay-wave dispersion relation, $\omega^2 = \Wsplay^2\left(\bk;\bnz\right)$.
The phase $\Phi$ then satisfies the
linearized eikonal equation
\[
\Phi_t^2 - a^2\left(\nabla \Phi;\bnz\right) = 0.
\]
We assume that $\Phi(\bx,t)$ is single-valued and
caustics do not arise.

We find that (see Appendix \ref{app:splay})
\[
\bn_1(\theta,\bx,t) =  u(\theta,\bx,t)\, \Rsplay(\bx,t),
\]
where $\Rsplay(\bx,t)$ is defined  by \eq{Rsplay}. The scalar wave-amplitude
function $u(\theta,\bx,t)$ satisfies the  HS-equation
\be
\left(u_{t} + \ba\cdot\nabla u+ \Gamma u u_\theta
+ P u\right)_\theta =
\frac{1}{2} \Gamma u_\theta^2.
\label{hseq}
\ee
In this equation, $\ba(\bx,t)$ is the linearized group velocity vector ($\ba=\nabla_\bk \omega$),
\[
\ba = \frac{1}{\omega}\left[\alpha \bk - (\alpha-\gamma) \left(\bk\cdot \bnz\right) \bnz\right] ,
\]
$\Gamma(\bx,t)$ is the genuine-nonlinearity coefficient ($\Gamma = \nabla_{\bnz} \omega\cdot\Rsplay$),
\[
\Gamma = -\left(\frac{\alpha - \gamma}{\omega}\right)\left(\bk\cdot\bnz\right)
\left[ k^2 - \left(\bk\cdot \bnz\right)^2\right],
\]
and $P(\bx,t)$ is given by
\[
P = \frac{\left(\omega R^2\right)_t + \div\left(\omega R^2 \ba\right)}{2\omega R^2},
\]
where $R^2 = k^2 - \left(\bk\cdot \bnz\right)^2$.
Equation \eq{hseq} follows from the variational principle
\[
\delta \int \frac{1}{2}\omega R^2 \left[u_\theta\left(u_t + \ba\cdot \nabla u\right)
+ \Gamma u u_\theta^2\right]\,d\theta d\bx dt = 0.
\]

Introducing a derivative along the rays associated with $\Phi$,
\[
\partial_s = \partial_t + \ba\cdot\nabla,
\]
we may write \eq{hseq} as an evolution equation for $u$ along a ray,
\be
\left(u_{s} + \Gamma u u_\theta
+ P u\right)_\theta =
\frac{1}{2} \Gamma u_\theta^2.
\label{rayevsplay}
\ee
If $J(\bx,t)$ is a non-zero ray-density function such that
$J_t + \div \left(J\ba\right) = 0$,
then we have
\[
P  = \frac{\left(\omega R^2\right)_s}{2\omega R^2} - \frac{J_s}{2 J},
\]
and the change of variables
\[
u\mapsto \sqrt{\frac{\omega R^2}{J}} \,u,
\qquad
\partial_\theta\mapsto\partial_x,
\qquad
\partial_s \mapsto \sqrt{\frac{J}{\omega R^2}}\Gamma\, \partial_t
\]
reduces \eq{rayevsplay} to \eq{HSeq}.

\subsection{Twist waves}

We consider the propagation of weakly nonlinear twist waves with
non-zero wavenumber
vector $\bk$ through an unperturbed director field $\bnz$. We assume that
$\alpha\ne \beta$ and $\bk$ is not parallel to $\bn_0$. These conditions
ensure the strict hyperbolicity of the system.

As a result of their linear degeneracy, the effect of nonlinearity on the twist
waves is cubic.
We therefore look for an asymptotic solution of \eq{dirPDE} of the form
\bea
&&\bn^\eps(\bx,t) = \bn\left(\frac{\Phi(\bx,t)}{\eps},t;\eps\right),
\label{nexpeps}
\\
&&\bn(\theta,\bx,t;\eps) = \bnz(\bx,t) + \eps^{1/2} \bn_1(\theta,\bx,t)
+ \eps \bn_2(\theta,\bx,t) + O(\eps^{3/2})
\label{nexp}
\eea
as $\eps \to 0^+$, where $\bnz$ is a solution of \eq{dirPDE}. We assume
that the local frequency and wavenumber \eq{locfreq} satisfy the linearized
dispersion relation for twist waves,
$\omega^2 = \Wtwist^2\left(\bk;\bnz\right)$,
so that the phase $\Phi$ satisfies the eikonal equation
\[
\Phi_t^2 - \Wtwist^2\left(\nabla \Phi;\bnz\right) = 0.
\]

We find that (see Appendix \ref{app:twist})
\bea
\bn_1 &=& v \Rtwist,
\label{nsol}
\\
\bn_2 &=& - \left(\frac{\beta - \gamma}{\alpha-\beta}\right) \left(\bk\cdot\bnz\right) u \Rsplay
+ \frac{1}{2}v^2 \left(\bk \times \Rtwist\right),
\label{n2sol}
\eea
where $\Rsplay(\bx,t)$, $\Rtwist(x,t)$ are defined in \eq{Rsplay}, \eq{Rtwist} and
the scalar amplitude-functions
$u(\theta,\bx,t)$, $v(\theta,\bx,t)$ satisfy
\bea
&&\left( v_t +\bb\cdot\nabla v + \Lambda u v_\theta + Q v\right)_\theta = 0,
\label{ueq}
\\
&&u_{\theta\theta} = v_\theta^2.
\label{veq}
\eea
Here, $\bb(\bx,t)$ is the group velocity vector,
\[
\bb = \frac{1}{\omega}\left[\beta \bk - (\beta-\gamma) \left(\bk\cdot \bnz\right) \bnz\right],
\]
$\Lambda(\bx,t)$ is given by
\[
\Lambda = \frac{(\beta-\gamma)^2}{\omega(\alpha-\beta)} \left(\bk\cdot\bnz\right)^2
\left[k^2 - \left(\bk\cdot\bnz\right)^2\right],
\]
and $Q(\bx,t)$ is given by
\[
Q = \frac{\left(\omega S^2\right)_t + \div\left(\omega S^2 \bb\right)}{2\omega S^2},
\]
where $S^2 = k^2 - \left(\bk\cdot \bnz\right)^2$. (With the normalization we adopt for
$\Rsplay$ and $\Rtwist$, we have $R^2 = S^2$.)
Equations \eq{ueq}--\eq{veq} follow from the variational principle
\[
\delta \int \frac{1}{2}\omega S^2 \left[v_\theta\left(v_t + \bb\cdot \nabla v\right)
+ \Lambda\left( u v_\theta^2+\frac{1}{2}u_\theta^2\right)\right]\,d\theta d\bx dt = 0.
\]

The solution \eq{nexpeps}--\eq{n2sol} consists of a leading-order
twist wave with amplitude $v$ and a higher-order forced splay wave with amplitude $u$
that propagates at the twist-wave velocity. The amplitude and frequency
of the forced splay wave are $O(\eps)$ and $O(1/\eps)$, respectively,
which is the same scaling as in the weakly nonlinear
solution for a free splay wave given in Section~\ref{sec:splay}.
The term proportional to $v^2$ in $\bn_2$ ensures that $\bn$
is a unit vector up to the first order in $\eps$.

The coefficient $\Lambda$ of the nonlinear
term in \eq{ueq} is non-zero if $\beta \ne \gamma$ and $\bk$ is
not parallel or orthogonal to $\bnz$.
If $\bk$, $\bnz$ are constant and $\bk$ is orthogonal to $\bnz$, then there are
exact large-amplitude traveling twist-wave solutions \cite{ericksen,ericksen1,shahinpoor}
and no weakly nonlinear effects arise.
If $\bk$ is parallel to $\bnz$, then there is a loss of strict
hyperbolicity and one obtains a system of asymptotic equations instead of a scalar equation.
We consider this case in Section~\ref{sec:rot}.

Introducing a derivative along the rays associated with $\Phi$,
\[
\partial_s = \partial_t + \bb\cdot\nabla,
\]
we may write \eq{ueq}--\eq{veq} as
\be
\left( v_s + \Lambda u v_\theta + Q v\right)_\theta = 0,
\qquad u_{\theta\theta} = v_\theta^2.
\label{rayevtwist}
\ee
If $K(\bx,t)$ is a non-zero ray-density function such that
$K_t + \div \left(K\bb\right) = 0$,
then we have
\[
Q  = \frac{\left(\omega S^2\right)_s}{2\omega S^2} - \frac{K_s}{2 K},
\]
and the change of variables
\[
u\mapsto\frac{\omega S^2}{K} \,u,
\quad
v\mapsto \sqrt{\frac{\omega S^2}{K}} \,v,
\qquad
\partial_\theta\mapsto\partial_x,
\qquad
\partial_s \mapsto \frac{K\Lambda}{\omega S^2}\, \partial_t
\]
reduces \eq{rayevtwist} to \eq{uveq}.

\subsection{Polarized waves}
\label{sec:rot}

The equations of motion \eq{dirPDE} are
invariant under spatial rotations and reflections that leave $\bn$ fixed.
As a consequence of this invariance, there is
a loss of strict hyperbolicity and genuine nonlinearity
for waves that propagate in the same direction as $\bn$.
The resulting polarized orientation waves are described by a cubically
nonlinear, rotationally invariant asymptotic equation.
An analogous phenomenon occurs for rotationally invariant waves
in first-order hyperbolic systems of conservation laws \cite{bh}.

We suppose that the unperturbed director field
$\bnz$ is constant, and consider waves with constant
wavenumber vector $\bk = k\bnz$ parallel to $\bnz$.
We look for an asymptotic solution $\bn^\eps$ of \eq{dirPDE} of the form
\beaa
&&\bn^\eps = \bn\left(\frac{\bk\cdot\bx - \omega t}{\eps},\bx,t;\eps\right),
\\
&&\bn(\theta,\bx,t;\eps) = \bnz + \eps^{1/2} \bn_1(\theta,\bx,t) + \eps \bn_2(\theta,\bx,t)
+ O(\eps^{3/2})
\eeaa
as $\eps \to 0^+$.

We find that  $\omega^2 = \gamma k^2$ (see Appendix \ref{app:pol}), and
\beaa
&&\bn_1 = \bu
\qquad\qquad\mbox{where $\bnz\cdot \bu = 0$},
\\
&&\bn_2 = -\frac{1}{2} \left(\bu\cdot\bu\right) \bnz.
\eeaa
The leading-order perturbation $\bu(\theta,\bx,t)$ satisfies the equation
\bea
&&\bu_{\theta t} + \frac{\omega}{k}\bnz\cdot\nabla \bu_\theta
+ \frac{(\alpha-\beta)k^2}{2\omega}\left( \bu\cdot\bu_\theta\right)_{\theta} \bu
\nonumber
\\
&&\qquad\qquad\quad+ \frac{(\beta-\gamma)k^2}{2\omega}
\biggl\{\left[\left(\bu\cdot\bu\right) \bu_{\theta}\right]_\theta - \left(\bu_\theta\cdot\bu_\theta\right) \bu\biggr\}
= 0,
\label{poleq}
\eea
which is derived from the variational principle
\beaa
&&\delta \int \frac{1}{2}\left\{\bu_\theta\cdot\left(\bu_t + \frac{\omega}{k}\bnz\cdot\nabla\bu\right)
+  \frac{(\alpha-\beta) k^2}{2\omega}\left(\bu\cdot\bu_\theta\right)^2\right.
\\
&&\qquad\qquad\qquad
\left.+  \frac{(\beta-\gamma) k^2}{2\omega}(\bu\cdot\bu) \left(\bu_\theta\cdot\bu_\theta\right)
\right\}\,d\theta d\bx dt = 0.
\eeaa

Making the change of variables
\[
\partial_t + \frac{\omega}{k}\bnz\cdot\nabla \mapsto \partial_t,\qquad \partial_\theta\mapsto \partial_x,
\]
and rescaling $\bu$, we can write  \eq{poleq} in a normalized form
for $\bu(x,t) \in \Rl^2$ as
\bea
&&\bu_{x t}
+ \left(\mu-\nu\right) \left(\bu\cdot\bu_x\right)_{x} \bu
+ \nu \left[\left(\bu\cdot\bu\right) \bu_{x}\right]_x
- \nu\left(\bu_x\cdot\bu_x\right) \bu
= 0,\qquad
\label{normpoleq}
\eea
where
\[
\mu =\frac{\alpha-\gamma}{\gamma},\qquad \nu = \frac{\beta-\gamma}{\gamma}.
\]
The corresponding variational principle is
\[
\delta \int \frac{1}{2}\left\{\bu_x\cdot\bu_t
+  (\mu-\nu)\left(\bu\cdot\bu_x\right)^2
+ \nu  (\bu\cdot\bu)\left(\bu_x\cdot\bu_x\right)
\right\}\,dx dt = 0.
\]

Writing $\bu = (u\cos v,u\sin v)$, we find that this variational principle
becomes
\be
\delta\int \frac{1}{2}\left\{ u_x u_t + \mu u^2 u_x^2 + u^2 \left(v_x v_t + \nu u^2 v_x^2\right)
\right\}\,dx dt = 0.
\label{redvar}
\ee
This result is consistent with what we obtain by expanding the one-dimensional variational
principle \eq{onedvar}--\eq{dircfs} as $\phi\to 0$, when
\[
a^2 \sim a_0^2 + (\alpha-\gamma)\phi^2, \quad b^2 \sim a_0^2 + (\beta-\gamma)\phi^2,
\quad q^2 \sim \phi^2,
\]
and making a unidirectional approximation $\partial_t \sim -a_0\partial_x$ in the
resulting Lagrangian.

The Euler-Lagrange equations for \eq{redvar} are
\bea
&&\left(u_t + \mu u^2u_x\right)_x - \mu u u_x^2 - u v_x\left(v_t+2\nu u^2v_x\right) = 0,
\label{pol1}
\\
&&\left(v_t + \nu u^2v_x\right)_x + 2\nu uu_x v_x + \frac{1}{u}\left(u_x v_t + u_t v_x\right) = 0.
\label{pol2}
\eea
This system is a coupled pair of wave equations for
$u$ and $v$. The radial mode $u$ has velocity $\mu u^2$, so it is genuinely
nonlinear when $u\ne0$, while the angular mode $v$ has velocity $\nu u^2$,
so it is linearly degenerate.
If $v$ is constant, corresponding to a plane-polarized wave, we recover the scalar cubic equation
\eq{cubic1} for $u$. Nonlinear circularly polarized waves do not exist, however, since variations in the
angular variable $v$ force variations in the radial variable $u$.


\section{One-dimensional equations}
\label{sec:oned}
\setcounter{equation}{0}

We consider a director field
\[
\bn(x,t) = \left(n_1(x,t), n_2(x,t), n_3(x,t)\right)
\]
that depends upon a single space variable $x = x_1$. Writing
\[
\bn = \left(\cos\phi, \sin\phi\cos\psi,\sin\phi\sin\psi\right),
\]
where $\phi(x,t)$, $\psi(x,t)$ are spherical polar angles,
we find that the variational principle \eq{varprinc}
becomes
\be
\delta \int \frac{1}{2}\left\{\phi_t^2 - a^2(\phi)\phi_x^2
 + q^2(\phi)\left[\psi_t^2 - b^2(\phi)\psi_x^2\right]\right\} \, dx dt = 0,
\label{onedvar}
\ee
with
\bea
&&a^2(\phi) =\alpha \sin^2\phi + \gamma\cos^2 \phi,
\nonumber\\
&&b^2(\phi) = \beta \sin^2\phi + \gamma\cos^2 \phi,
\label{dircfs}
\\
&&q^2(\phi) = \sin^2 \phi.
\nonumber
\eea

The Euler-Lagrange equation associated with \eq{onedvar} is
a system of wave equations,
\bea
&&\phi_{tt} - \left(a^2\phi_{x}\right)_x + aa^\prime\phi_x^2 + q^2 bb^\prime \psi_x^2
- qq^\prime\left(\psi_t^2 - b^2\psi_x^2\right) = 0,
\label{wavea}\\
&&
\psi_{tt} - \left(b^2 \psi_{x}\right)_x
+\frac{2q^\prime}{q}\left(\phi_t\psi_t - b^2\phi_x\psi_x\right) = 0,
\label{waveb}
\eea
where the prime denotes a derivative with respect to $\phi$.
The angle $\phi$ corresponds
to splay waves and the angle $\psi$ to twist waves.
Both wave-speeds $(a,b)$ are functions only of $\phi$, and
the wave equation for $\phi$ is forced by terms that are proportional to
quadratic functions of derivatives of $\psi$.

\subsection{The initial value problem}

Hunter and Saxton \cite{hs} construct an
asymptotic solution $\phi(x,t;\eps)$ of the initial value problem for
the scalar wave equation \eq{scalarwaveeq}
in $-\infty < x < \infty$ and $t > 0$
with weakly nonlinear splay-wave initial data,
\beaa
\phi(x,0;\eps) = \phi_0 + \eps f\left(\frac{x}{\eps}\right),\qquad \phi_t(x,0;\eps) = g\left(\frac{x}{\eps}\right),
\eeaa
where $f$, $g$ have compact support, and $\eps$
is a small parameter.
The solution consists of a superposition of right and left moving weakly nonlinear
splay waves that originate from $x=0$, whose width in $x$ is of the order $\eps$. The splay
waves are separated by a slowly-varying,
small-amplitude perturbation
of the constant state $\phi_0$ that satisfies a linearized wave equation.

Here, we construct an asymptotic solution for $\phi(x,t;\eps)$, $\psi(x,t;\eps)$
of the system \eq{wavea}--\eq{waveb} with weakly nonlinear twist-wave initial data,
\bea
&&\phi(x,0;\eps) = \phi_0,\qquad\qquad\qquad\qquad \phi_t(x,0;\eps) = 0,
\nonumber
\\
&&\psi(x,0;\eps) = \eps^{1/2}f\left(\frac{x}{\eps}\right),
\qquad \psi_t(x,0;\eps) = \eps^{-1/2}g\left(\frac{x}{\eps}\right),
\label{twistic}
\eea
where $f$ and $g$ are smooth, compactly supported functions.
We suppose, without loss of generality,
that the unperturbed constant value of $\psi$ is equal to
zero.

We use a $0$-subscript on a function of $\phi$ to denote evaluation
at $\phi=\phi_0$. We assume that $a_0\ne b_0$,
so that \eq{wavea}--\eq{waveb} is strictly hyperbolic at $\phi_0$. The system
will then remain strictly hyperbolic, at least in some time-interval of the order one. We further
assume that $b_0^\prime \ne 0$, otherwise the leading-order nonlinear effects studied below vanish.
In the case of the director-field wave speeds \eq{dircfs},
these assumptions mean that $\phi_0 \ne n\pi/2$ for $n\in \Ir$.

Although $\phi$ is initially constant, it does not remain so.
As we will show, the weakly nonlinear twist wave, whose initial energy
\[
\frac{1}{2} \int_{-\infty}^\infty\left\{ \psi_t^2 + b^2(\phi) \psi_x^2\right\}\, dx
= \frac{1}{2} \int_{-\infty}^\infty\left\{g^2(\theta) + b_0^2f_\theta^2(\theta)\right\}\, d\theta
\]
is of the order one, generates a slowly-varying
`outer' splay-wave whose amplitude is of the order one.

We will construct an asymptotic solution of this initial value problem as follows.
\begin{enumerate}
\item In a short initial layer, when $t=O(\eps)$, we use linearized theory. The initial data
splits up into right and left moving twist waves.
\item For $t=O(1)$, nonlinear
effects become important, and we use the method of matched asymptotic expansions,
with different expansions for the twist and splay waves.
\begin{enumerate}
\item The twist waves are small-amplitude, localized waves,
which vary on a spatial scale of the order $\eps$.
We describe them by means of the
weakly nonlinear asymptotic equations for twist waves derived above. We call these the `inner' solutions.
\item Away from the twist waves, the leading-order solution is a large-amplitude
splay wave, which varies on a spatial scale of the order 1. We call this the `outer' solution.
\item We obtain jump conditions for the `outer' splay-wave solution
across the trajectories of the right and left moving twist waves
by matching it with the `inner' twist-wave
solutions.
\item We obtain initial data for the nonlinear solution
by matching it with the linearized solution as $t\to 0^+$.
\end{enumerate}
\end{enumerate}

\begin{figure}
\centerline{
\epsfxsize=2.5in
\epsfysize=2.5in
\epsfbox{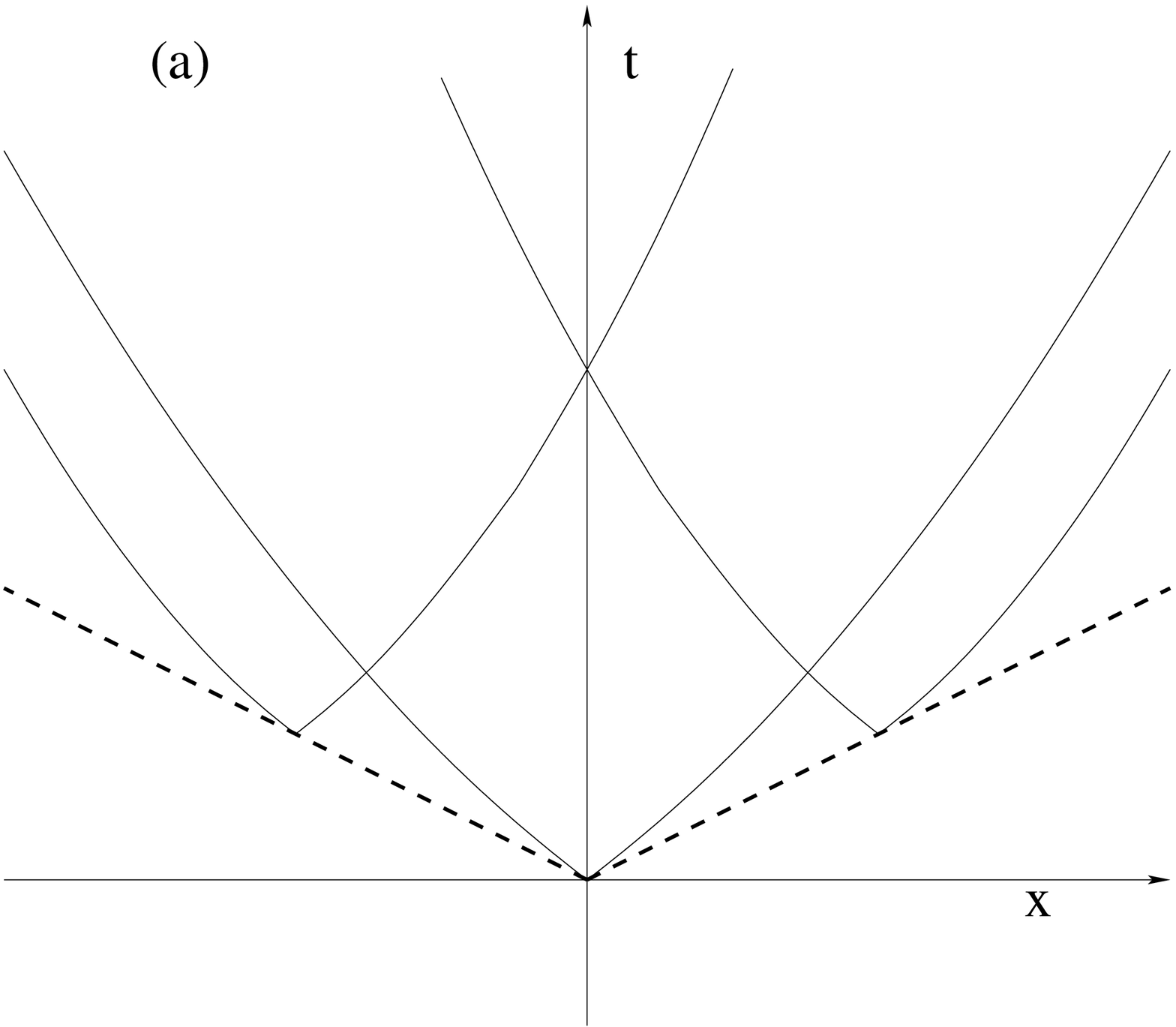}
\hskip0.2in
\epsfxsize=2.5in
\epsfysize=2.5in
\epsfbox{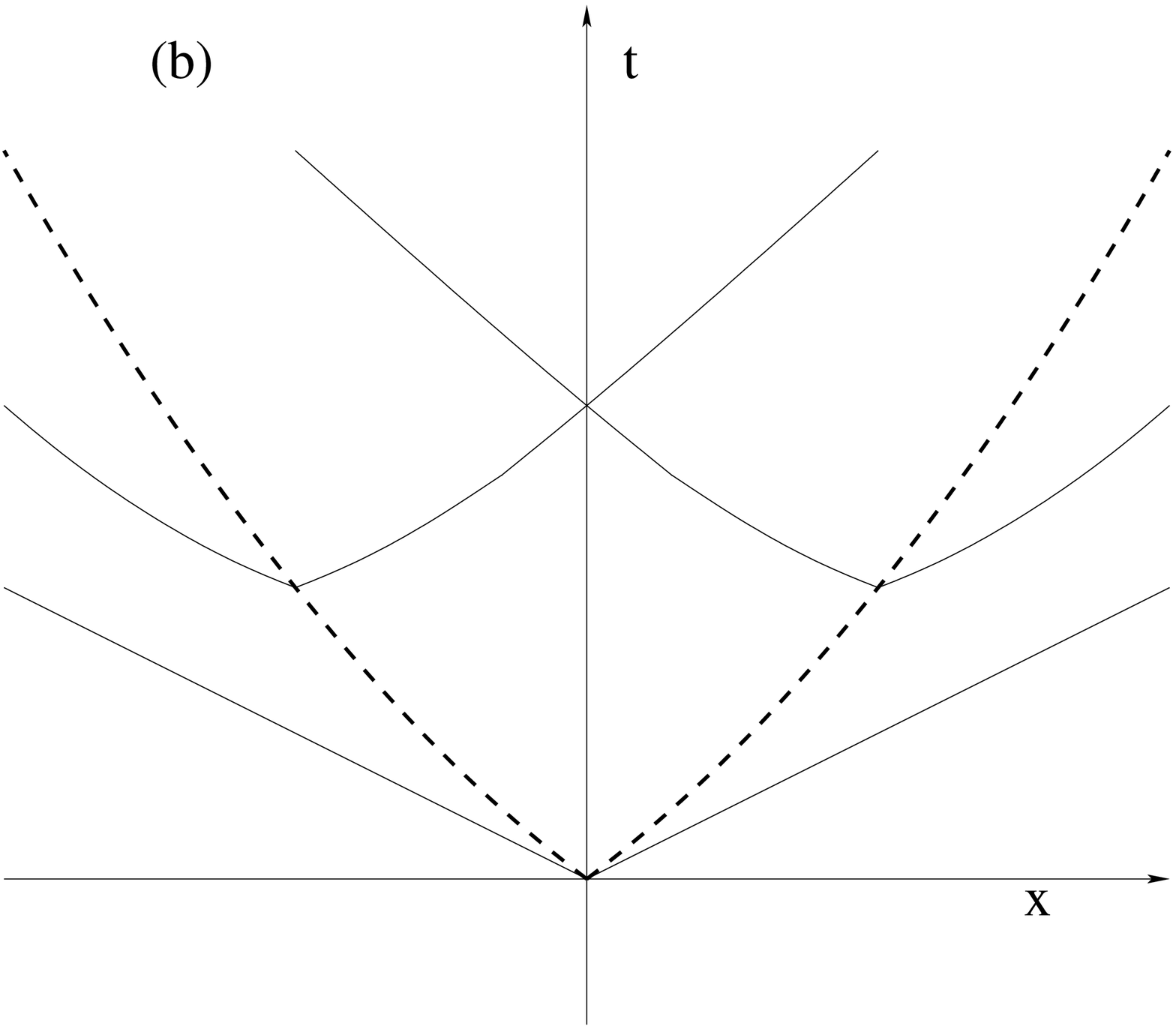}}


\caption{Characteristic structure for the solution of the initial value problem
\eq{wavea}--\eq{twistic}:
(a) Fast twist waves; (b) Slow twist waves.
The dashed lines are the trajectories of the twist waves, and the solid lines
are the characteristic curves associated with the splay wave. \label{charivp}}
\end{figure}

The matching between the inner twist waves and the outer splay wave,
and the structure of the nonlinear solution,
depend on whether the twist waves are faster
or slower than the splay waves (see Figure~\ref{charivp}).
In deriving the jump conditions for the splay wave
below, we will consider these cases separately.
For liquid crystals, we have $\beta < \alpha$, meaning that
twist deformations are not as `stiff' as splay deformations,
and the twist waves are slower.

\subsection{Initial layer}

First, we consider the solution of \eq{wavea}--\eq{waveb}, \eq{twistic} in a short initial layer
when $t=O(\eps)$.
Because of the finite propagation speed of the system, the solution is constant,
with $\phi=\phi_0$, $\psi=0$,
outside an interval of width of the order $\eps$ containing $x=0$.

Near $x=0$, we look for an asymptotic solution of the form
\[
\phi = \Phi_0(X,T) + O(\eps),\quad \psi = \eps^{1/2} \Phi_1(X,T) + O(\eps),
\quad
X = \frac{x}{\eps},\quad T = \frac{t}{\eps}.
\]
Using this expansion in \eq{wavea}--\eq{twistic},
we find that $\Phi_0$, $\Psi_1$ satisfy
\beaa
&&\Phi_{0TT} - \left(a^2(\Phi_0) \Phi_{0X}\right)_X  + a(\Phi_0) a^\prime(\Phi_0) \Phi_{0X}^2= 0,
\\
&&\Psi_{1TT} - \left(b^2(\Phi_0) \Psi_{1X}\right)_X
+\frac{2q^\prime(\Phi_0)}{q(\Phi_0)}\left(\Phi_{0T}\Psi_{1T}
- b^2(\Phi_0)\Phi_{0X}\Psi_{1X}\right) = 0,
\\
&&\Phi_0(X,0) = \phi_0,\qquad\qquad \Phi_{0T}(X,0) = 0,
\\
&&\Psi_1(X,0) = f\left(X\right),
\qquad \Psi_{1T}(X,0) = g\left(X\right).
\eeaa
The initial value problem for $\Phi_0$ has the constant solution
$\Phi_0 = \phi_0$, and therefore $\Psi_1$ satisfies the linear wave equation
\beaa
&&\Psi_{1TT} - b_0^2 \Psi_{1XX} = 0,
\\
&&\Psi_1(X,0) = f\left(X\right),
\quad
\Psi_{1T}(X,0) = g\left(X\right).
\eeaa
The solution is
\bea
&&\Psi_1(X,T) = F_R\left(X-b_0 T\right) + F_L\left(X+b_0 T\right),
\label{dalembert}
\\
&&F_R(\theta) = \frac{1}{2} f(\theta) + \frac{1}{2 b_0} \int^\infty_\theta g(\xi)\, d\xi,
\label{defF}
\\
&&F_L(\theta) = \frac{1}{2} f(\theta) - \frac{1}{2 b_0} \int^\infty_\theta g(\xi)\, d\xi.
\label{defG}
\eea
Since $f$, $g$ have compact support, the functions $F_{R\theta}$, $F_{L\theta}$ have compact support.

\subsection{Twist waves}

Next, we consider the propagation of the twist-waves for $t=O(1)$
through a possibly non-uniform splay-wave field $\phi(x,t)$.
For definiteness, we consider the right-moving twist wave, which
moves with velocity $b(\phi)$.
The trajectory $x=s_R(t)$ of the wave satisfies
\[
\frac{d s_R}{dt} = b_R,
\]
where an $R$-subscript on a function of $\phi$
denotes evaluation at $\phi=\phi_R(t)$, with
$\phi_R(t) = \phi\left(s_R(t),t\right)$. For the
initial value problem, we have
$s_R(0) = 0$ and $\phi_R(0) = \phi_0$.

We introduce a stretched inner variable
near this trajectory,
\[
\theta = \frac{x - s_{R}(t)}{\eps}.
\]
We find that the weakly nonlinear
twist-wave solution of \eq{wavea}--\eq{waveb} is
\bea
&&\phi = \phi_R(t) + \eps \phi_2(\theta,t)+ O(\eps^{3/2}),
\label{phisol}
\\
&&\psi = \eps^{1/2} \psi_1(\theta,t) + O(\eps),
\label{psisol}
\eea
where $\phi_2$, $\psi_1$ satisfy
\bea
&&\left[\psi_{1t} + b_R^\prime \phi_2 \psi_{1\theta} + \left(\frac{b_{Rt}}{2b_R}
+ \frac{q_{Rt}}{q_R}\right)\psi_1\right]_\theta = 0,
\label{ibvp1}
\\
&&\phi_{2\theta\theta} = \left(\frac{q_R^2 b_R b_R^\prime}{a_R^2 - b_R^2}\right) \psi_{1\theta}^2.
\label{ibvp2}
\eea
We may transform \eq{ibvp1}--\eq{ibvp2} into \eq{uveq} by
a suitable change of variables, in which $\theta$ corresponds to $x$, $\psi_1$ to $v$,
and $\phi_2$ to $u$.

Matching the twist-wave solution \eq{psisol} as $t\to 0^+$
with the linearized solution \eq{dalembert}
as $T \to \infty$, we get the initial condition
\be
\psi_1(\theta,0) = F_R(\theta),
\label{ic3}
\ee
where $F_R$ is given in \eq{defF}.
Equations \eq{ibvp1}--\eq{ibvp2} are supplemented with
suitable boundary conditions for $\phi_2$ and $\psi_1$ at $\theta = \pm\infty$,
which we consider further below.

The main result we need in order to obtain equations for the `outer' splay
wave solution is the following jump condition for $\phi_{2\theta}$ across the twist wave.
We give a complete solution of the twist-wave equations in Section~\ref{sec:char}.

\begin{proposition}
Suppose that $\phi_2$, $\psi_1$ are smooth solutions of \eq{ibvp1}--\eq{ibvp2}
such that $\psi_{1\theta}(\cdot,t)$ has compact support. Then
\be
\frac{d}{dt}\left\{ \left(\frac{a_R^2 - b_R^2}{b_R^\prime}\right) \left[\phi_{2\theta}\right]\right\}
+ \left(\frac{a_R^2 - b_R^2}{2}\right) \left[\phi_{2\theta}^2\right] = 0,
\label{jumpcon}
\ee
where $[\cdot]$ denotes the jump from $\theta=-\infty$ to $\theta=\infty$.
\end{proposition}

\medskip\noindent
\textbf{Proof.}
Since $\psi_{1\theta}$ has compact support, it follows from \eq{ibvp2} that $\phi_2$ is a linear
function of $\theta$ for large negative and positive values of $\theta$. Moreover,
\be
\left[\phi_{2\theta}\right] = \left(\frac{q_R^2 b_R b_R^\prime}{a_R^2 - b_R^2}\right)
\int_{-\infty}^\infty \psi_{1\theta}^2\, d\theta.
\label{intjump}
\ee
Multiplying \eq{ibvp1} by $\psi_{1\theta}$ and using \eq{ibvp2} to rewrite the result,
we get
\[
\left(\psi_{1\theta}^2\right)_t
+ \left[b_R^\prime \phi_2 \psi_{1\theta}^2 +
\frac{1}{2}\left(\frac{a_R^2 - b_R^2}{q_R^2 b_R}\right) \phi_{2\theta}^2\right]_\theta
+\left(\frac{b_{Rt}}{b_R} + \frac{2q_{Rt}}{q_R}\right) \psi_{1\theta}^2 = 0.
\]
Integrating this equation with respect to $\theta$, and using the fact that $\psi_{1\theta}$
has compact support, we find that
\[
\frac{d}{dt} \int_{-\infty}^\infty\psi_{1\theta}^2 \, d\theta
+\frac{1}{2}\left(\frac{a_R^2 - b_R^2}{q_R^2 b_R}\right) \left[\phi_{2\theta}^2\right]
+ \left(\frac{b_{Rt}}{b_R} + \frac{2q_{Rt}}{q_R}\right) \int_{-\infty}^\infty \psi_{1\theta}^2\, d\theta = 0.
\]
Using \eq{intjump} to eliminate the integral of $\psi_{1\theta}^2$ from this equation,
and rewriting the result, we get \eq{jumpcon}. \qquad$\Box$

\medskip
For the coefficients in \eq{dircfs}, the jump condition \eq{jumpcon} becomes
\[
\frac{d}{dt}\left\{ b_R\tan \phi_{R} \left[\phi_{2\theta}\right]\right\}
+ \frac{1}{2} \left(\beta-\gamma\right) \sin^2\phi_{R} \left[\phi_{2\theta}^2\right] = 0.
\]

\subsection{Matching: fast twist waves}
\label{sec:fast_twist}

This case corresponds to $0 < \alpha < \beta$, when $0 < a< b$.
Since the twist waves are faster than the splay waves,
they propagate into a constant state $\phi=\phi_0$, $\psi=0$
ahead of them,
and generate splay waves behind them. (See Figure~\ref{charivp}(a).)
It follows that the right and left moving
twist waves move at a constant velocity along the trajectories $x=b_0t$ and $x=-b_0t$, respectively.

We consider the right-moving twist wave for definiteness. The appropriate inner variable
is then
\be
\theta = \frac{x-b_0t}{\eps}.
\label{deftheta}
\ee
The weakly nonlinear solution inside the twist wave must match for
large positive $\theta$ with
the constant initial state ahead of the wave. This
condition implies that $\phi_2(\theta,t)$ and $\psi_1(\theta,t)$
in \eq{phisol}--\eq{psisol},
with $\phi_R = \phi_0$, satisfy
the boundary conditions
\be
\phi_2(\infty,t) = 0,
\qquad
\phi_{2\theta}(\infty, t) = 0,
\qquad
\psi_1(\infty,t) = 0.
\label{icbc}
\ee

The derivative $\phi_{2\theta}$ jumps from zero at $\theta=\infty$
to a value
\[
\phi_{2\theta}(-\infty,t) = \sigma_R(t)
\]
at $\theta = -\infty$. It
follows from the jump condition \eq{jumpcon}, with $a_R=a_0$, $b_R = b_0$, $b_R^\prime=b_0^\prime$
constants, that $\sigma_R$ satisfies
\[
\frac{d\sigma_R}{dt}
+ \frac{1}{2} b_0^\prime\sigma_R^2 = 0.
\]
From \eq{ibvp2} and \eq{ic3}, we have
$\sigma_R(0) = \sigma_{R0}$ where
\[
\sigma_{R0} = -\left(\frac{q_0^2 b_0 b_0^\prime}{a_0^2 - b_0^2}\right) \int_{-\infty}^\infty F_{R\theta}^2\, d\theta .
\]
The solution of this Riccati equation,
\be
\sigma_R(t) = \frac{\sigma_{R0}}{1 + \sigma_{R0} b_0^\prime t/2},
\label{solsig}
\ee
is defined for all $t \ge 0$, since $\sigma_{R0} b_0^\prime \ge 0$
when $0 < a_0 < b_0$.

In Section~\ref{sec:char}, we prove that when $a_0 < b_0$, equations
\eq{ibvp1}--\eq{ic3},  \eq{icbc}, with $a_R = a_0$ and so on,
have a smooth solution defined for all $t \ge 0$.
We note that the derivative $\phi_{2\theta}(-\infty,t)$ decays as $t\to \infty$. This
is a result of the fact that
the twist wave radiates energy away from it in the form of splay waves.
It also follows from
the solution that the twist wave is a rarefaction, in the sense that
its characteristics spread out with increasing time.

A similar analysis applies to the left-moving twist wave, in which
\[
\theta = \frac{x + b_0 t}{\eps},
\]
and $\phi_{2\theta}(\theta,t) = 0$ for $\theta$ sufficiently large and negative.
We find that $\phi_{2\theta} = \sigma_L$ for $\theta$
sufficiently large and positive, where
\be
\sigma_L(t) = \frac{\sigma_{L0}}{1-\sigma_{L0} b_0^\prime t/2},
\label{solmu}
\ee
with
\[
\sigma_{L0} = \left(\frac{q_0^2 b_0 b_0^\prime}{a_0^2 - b_0^2}\right)
\int_{-\infty}^\infty F_{L\theta}^2 \, d\theta.
\]

These `inner' twist-wave solutions provide matching conditions for an `outer' splay-wave solution
$\phi(x,t)$.
Using \eq{phisol} with $\phi_R = \phi_0$ to rewrite the condition for the right-moving twist wave,
\[
\phi_{2\theta}(\theta,t) \sim \sigma_R(t)\qquad \mbox{as $\theta\to -\infty$}
\]
where $\theta$ is given by \eq{deftheta}, in terms of the outer solution $\phi(x,t)$, and equating
the outer limit of the inner solution with the inner limit of the outer solution,
we find that
\[
\phi_x(x,t) \sim \sigma_R(t) \qquad\mbox{as $x\to b_0t^-$}.
\]
Furthermore, the leading-order outer solution for $\phi$ is continuous across the twist-wave,
and $\psi$ is higher-order
in $\eps$. We obtain a condition for $\phi$ as $x\to-b_0t^+$ in an analogous way.

Summarizing these results for the leading-order outer
splay-wave solution $\phi(x,t)$,
we find that $\phi=\phi_0$ is constant if $x > b_0t$
or $x < -b_0t$. Inside the region $-b_0 t < x < b_0 t$,
we find that
$\phi$ satisfies \eq{scalarwaveeq}, with data
on the space-like lines $x = \pm b_0t$ given by
\beaa
&&\phi(b_0t, t)  = \phi_0, \qquad\quad \phi_x(b_0t,t) = \sigma_R(t),
\\
&&\phi(-b_0t, t)  = \phi_0, \qquad \phi_x(-b_0t,t) = \sigma_L(t),
\eeaa
where $\sigma_R$, $\sigma_L$ are given by \eq{solsig}, \eq{solmu}, respectively.

Since the initial value problem
for \eq{scalarwaveeq} is well-posed,
this Cauchy problem is presumably solvable. The solution $\phi$ may form singularities,
in which case it would have to be continued by a weak solution.

\subsection{Matching: slow twist waves}

This case corresponds to $\alpha > \beta > 0$, when $a> b > 0$.
Since the twist waves are slower than the splay waves, they generate
splay waves both in front and behind them. As a result,
the twist waves are embedded inside a splay-wave
field. (See Figure \ref{charivp}(b).) The speeds of the twist waves depend on the splay-wave field,
leading to a free-boundary problem for the trajectories
of the twist waves, coupled with a wave
equation for the splay wave that is subject to
jump conditions across the twist-wave trajectories.

We will not write out detailed asymptotic equations for the weakly nonlinear twist waves
in this case, but we summarize the equations satisfied by the leading-order outer splay-wave solution
$\phi(x,t)$. The main point is the derivation of jump conditions
for $\phi_x$ across the twist-wave trajectories.

The right and left moving twist waves are located at $x = s_{R}(t)$ and $x=s_{L}(t)$, respectively, where
\bea
&&\frac{d s_R}{dt}(t) = b\left(\phi\left( s_R(t), t\right)\right),\quad
\frac{d s_L}{dt}(t) = -b\left(\phi\left( s_L(t), t\right)\right),
\label{fb2}
\\
&&s_R(0) = 0,\qquad\qquad\qquad\quad  s_L(0) = 0.
\label{fb3}
\eea
The solution $\phi$ is continuous across $x = s_{R}(t)$ and $x=s_L(t)$, so that
\be
\left[\phi\right]_R = 0,\qquad
\left[\phi\right]_L = 0.
\label{fb4}
\ee
Here, and below, we use
$[\cdot]_R$, $[\cdot]_L$ to denote the jumps across $x=s_R(t)$, $x=s_L(t)$, respectively,
meaning that
\beaa
\left[\phi\right]_R(t) &=& \lim_{x\to s_R(t)^+} \phi(x,t) -  \lim_{x\to s_R(t)^-} \phi(x,t),
\\
\left[\phi\right]_L(t) &=& \lim_{x\to s_L(t)^+} \phi(x,t) -  \lim_{x\to s_L(t)^-} \phi(x,t).
\eeaa

Considering the right-moving twist wave for definiteness, we have
\bea
&&\phi_{2\theta}(\theta,t) \to \sigma_+(t)\qquad \mbox{as $\theta\to \infty$},
\label{slowjc1}
\\
&&\phi_{2\theta}(\theta,t) \to \sigma_-(t)\qquad \mbox{as $\theta\to -\infty$},
\label{slowjc2}
\eea
for some functions $\sigma_+(t)$, $\sigma_-(t)$. The corresponding matching conditions
for $\phi(x,t)$ are
\be
\lim_{x\to s_R^+(t)} \phi_x(x,t) = \sigma_+(t),\qquad
\lim_{x\to s_R^-(t)} \phi_x(x,t) = \sigma_-(t).
\label{matchjc}
\ee
From \eq{jumpcon}, \eq{slowjc1}--\eq{slowjc2}, \eq{matchjc},
and the analogous equations for the left-moving twist wave, we find that
$\phi_x$ satisfies the following jump conditions
across the twist waves:
\bea
&&\frac{d}{dt}\left\{ \left(\frac{a_R^2 - b_R^2}{b_R^\prime}\right) \left[\phi_{x}\right]_R\right\}
+ \left(\frac{a_R^2 - b_R^2}{2}\right) \left[\phi_{x}^2\right]_R = 0,
\label{fb7}
\\
&&
\frac{d}{dt}\left\{ \left(\frac{a_L^2 - b_L^2}{b_L^\prime}\right) \left[\phi_{x}\right]_L\right\}
-\left( \frac{a_L^2 - b_L^2}{2}\right) \left[\phi_{x}^2\right]_L = 0.
\label{fb8}
\eea
Furthermore, from \eq{ic3} and \eq{intjump}, and their analogs for left-moving waves,
we get the initial conditions
\bea
\left[\phi_x\right]_R(0) &=& \left(\frac{q_0^2 b_0 b_0^\prime}{a_0^2 - b_0^2}\right)
\int_{-\infty}^\infty F_{R\theta}^2\, d\theta,
\label{fb9}
\\
\left[\phi_x\right]_L(0) &=& \left(\frac{q_0^2 b_0 b_0^\prime}{a_0^2 - b_0^2}\right)
\int_{-\infty}^\infty F_{L\theta}^2\, d\theta.
\label{fb10}
\eea

Finally, matching the outer solution with the initial, linearized
solution, we find that $\phi(x,t)$ satisfies
the initial conditions
\be
\phi(x,0) = \phi_0,\qquad \phi_t(x,0) = 0, \qquad \mbox{for $x\ne 0$}.
\label{fb1}
\ee

Summarizing, we find that the free-boundary problem for $\phi$, $s_R$, $s_L$
consists of \eq{scalarwaveeq} for $\phi(x,t)$ in $-\infty < x < \infty$,
$t>0$ with the initial condition \eq{fb1}.
The functions $s_R(t)$, $s_L(t)$
satisfy \eq{fb2}--\eq{fb3}, and
$\phi(x,t)$ satisfies the jump conditions
\eq{fb4}, \eq{fb7}--\eq{fb10} across the curves
$x = s_R(t)$, $x = s_L(t)$.

We will not investigate this problem here.
We remark, however, that Proposition~\ref{prop:smoothsol} in Section~\ref{sec:char}
implies that there is a smooth
solution of the `inner' asymptotic equations for the weakly nonlinear twist wave
\eq{ibvp1}--\eq{ic3}, \eq{slowjc1}--\eq{slowjc2}
whenever the derivatives $\sigma_{\pm}(t) = \phi_x\left(s_{R}^\pm(t),t\right)$ of the `outer'
splay-wave solution of the free-boundary problem on either side of the
twist wave are smooth functions of time.

\section{Method of characteristics}
\label{sec:char}
\setcounter{equation}{0}

In this section, we solve \eq{uveq} by the method of characteristics. The
explicit nature of this solution is related to the complete
integrability of the equation, which is discussed in the next section.

\begin{proposition}
\label{prop:char_sol}
Let $(\xi,\tau)$ be characteristic coordinates for the PDE \eq{uveq}, where
$x=X(\xi,\tau)$, $t=\tau$, and write
$U(\xi,\tau) = u\left(X(\xi,\tau),\tau\right)$, $V(\xi,\tau) = v\left(X(\xi,\tau),\tau\right)$.
Then a formal solution of \eq{uveq} is given by
\be
U = X_\tau,\qquad V = F + G,\qquad X = - \int_0^\xi\frac{F_\xi^2(A+B)^2}{2A_\xi B_\tau}\,d\xi + H,
\label{char_sol}
\ee
where $A(\xi)$, $B(\tau)$, $F(\xi)$, $G(\tau)$, $H(\tau)$ are arbitrary functions.
\end{proposition}

\medskip\noindent
\textbf{Proof.}
Writing \eq{uveq} in terms of
characteristic coordinates $(\xi,\tau)$ in
which $\tau =t$ and $x_\tau = u$,
we find that the PDE becomes
\[
X_\tau = U,
\qquad
V_{\xi\tau} = 0,
\qquad
U_{\xi\xi} - \frac{J_{\xi}}{J} U_\xi = V_\xi^2,
\]
where $J(\xi,\tau)$ is the Jacobian $J = X_\xi$.
It follows that $V(\xi,\tau) = F(\xi) + G(\tau)$,
where $F$, $G$ are functions of integration, and
\[
J_\tau = U_\xi,
\qquad
U_{\xi\xi} - \frac{J_\xi}{J} U_\xi = F_\xi^2.
\]
The elimination of $U$ from these equations yields a PDE for $J$,
\[
J_{\xi\tau} - \frac{J_\xi J_\tau}{J} = F_\xi^2.
\]
Making the change of variables $\eta = \eta(\xi)$ where $\eta_\xi = F_\xi^2$, and
$J = -e^{-K}$, we find that this PDE transforms into an integrable Liouville equation,
\[
K_{\eta\tau} = e^{K}.
\]
The general solution is
\[
e^K = \frac{2 A_\eta B_\tau}{(A+B)^2},
\]
where $A(\eta)$ and $B(\tau)$ are arbitrary functions. Integrating the
equation $X_\xi = -e^{-K}$
with respect to $\xi$ , we find that $X$ is given by \eq{char_sol},
which proves the result. \qquad$\Box$


\medskip
Next, we consider \eq{uveq} in $-\infty < x < \infty$,
\bea
&&\left(v_t + u v_x\right)_x = 0,
\label{uveq1}
\\
&&u_{xx} = v_x^2,
\label{uveq2}
\eea
supplemented with the initial condition
\be
v(x,0) = F(x),
\label{uveqic}
\ee
and the boundary conditions
\be
u_x(\infty,t) = \sigma_+(t),\qquad u_x(-\infty,t) = \sigma_-(t).
\label{uxcond}
\ee

\begin{figure}
\centerline{
\epsfxsize=2.5in
\epsfysize=2.5in
\epsfbox{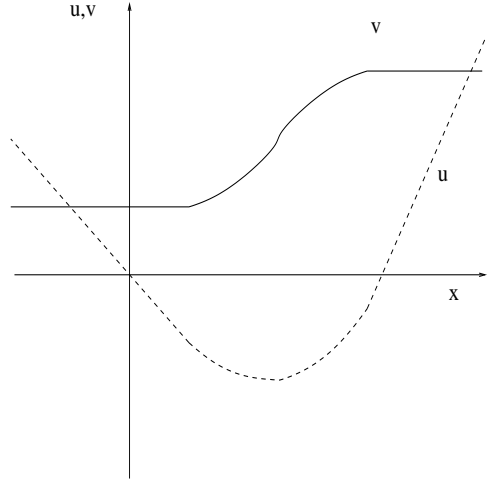}
}


\caption{Schematic structure of the solution of the IBVP \eq{uveq1}--\eq{uxcond}. \label{fig:soluveq}}
\end{figure}

We assume that $F$ is a smooth function and that $F_x$ has compact support.
Equation \eq{uveq1} then implies that $v_x(x,t)$ has compact
support in $x$ whenever a smooth solution exists,
and \eq{uveq2} implies that $u(x,t)$ is a linear function of $x$ for sufficiently
large positive and negative values of $x$. The boundary condition
\eq{uxcond} specifies the corresponding values of $u_x$.
We illustrate the structure of the solution
schematically in Figure~\ref{fig:soluveq}.

Computing the jump condition \eq{jumpcon} for \eq{uveq1}--\eq{uveq2}, we find that
\[
\frac{d}{dt} \left[u_x\right] + \frac{1}{2} \left[u_x^2\right] = 0,
\]
where $[\cdot]$ denotes the jump from $x=-\infty$ to $x=\infty$. It follows that
the data $\sigma_{\pm}$ must satisfy
\be
\frac{d\sigma_+}{dt} + \frac{1}{2} \sigma_+^2 = \frac{d\sigma_-}{dt} + \frac{1}{2} \sigma_-^2.
\label{jumpcons}
\ee
Moreover, integrating \eq{uveq2} with respect to $x$ at $t=0$ and using \eq{uveqic}--\eq{uxcond},
we find that
\be
\sigma_+(0) - \sigma_-(0) = \int_{-\infty}^\infty F_x^2 \, dx.
\label{sigzero}
\ee

The next proposition establishes the existence of smooth solutions of the IBVP
\eq{uveq1}--\eq{uxcond} for compatible data  $F$ and $\sigma_{\pm}$. The
solutions are not unique, since we may add an arbitrary
function of time to $v$, and an arbitrary
function of time to $u$ (together with
an appropriate time-dependent translation of the spatial coordinate $x$).
We can remove this non-uniqueness by specifying, for example,  $u$, $v$ as
functions of time at some value of $x$.

\begin{proposition}
\label{prop:smoothsol}
Suppose that $F:\Rl\to \Rl$ is a smooth function and that $F_x$ has compact support.
Also suppose that
$\sigma+, \sigma_-: [0,t_\ast) \to \Rl$
are smooth functions defined in some time
interval $0 \le t < t_\ast$, where $0 < t_\ast \le \infty$,
that satisfy \eq{jumpcons}, \eq{sigzero}.
Then there is a smooth solution of
\eq{uveq1}--\eq{uxcond}, defined in $-\infty < x < \infty$, $0\le t < t_\ast$.
\end{proposition}

\medskip\noindent
\textbf{Proof.}
If $F(x) = F_0$ is constant, then $\sigma_+(t) = \sigma_-(t)$, and a solution
is $u = x \sigma_+(t)$, $v = F_0$.
We therefore assume that $F$ is not constant.

We then
have from \eq{sigzero} that $\sigma_-(0) < \sigma_+(0)$. Equation
\eq{jumpcons} implies that
\[
\frac{d}{dt}\left(\sigma_+ - \sigma_-\right) + \frac{1}{2}\left(\sigma_+ + \sigma_-\right)
\left(\sigma_+ - \sigma_-\right) = 0,
\]
so $\sigma_-(t)< \sigma_+(t)$ for all $0  \le t < t_\ast$.

We choose constants $\eta_- < \eta_+ < 0$ such that
\[
\eta_+ - \eta_- = \int_{-\infty}^\infty F_\xi^2(\xi)\, d\xi,
\]
and define the function $\eta(\xi)$ by
\beaa
\eta(\xi) &=& \eta_- + \int_{-\infty}^\xi F_{\xi^\prime}^2(\xi^\prime)\, d\xi^\prime
\\
&=& \eta_+ - \int_\xi^{\infty} F_{\xi^\prime}^2(\xi^\prime)\, d\xi^\prime.
\eeaa
We also define a Jacobian $J(\xi,\tau)$ by
\be
J = \frac{\left[\left(\eta_+ - \eta\right)E_+ + \left(\eta - \eta_-\right)
E_-\right]^2}{\left(\eta_+ - \eta_-\right)\left(\sigma_+-\sigma_-\right)},
\label{defJ}
\ee
where
\bea
E_+(\tau) &=& \exp\left\{ -\frac{1}{4} \int_0^\tau\left[\sigma_+(\tau^\prime) - \sigma_-(\tau^\prime)\right]\, d\tau^\prime\right\},
\label{defE+}
\\
E_-(\tau) &=& \exp\left\{ +\frac{1}{4} \int_0^\tau\left[\sigma_+(\tau^\prime) - \sigma_-(\tau^\prime)\right]\, d\tau^\prime\right\}.
\label{defE-}
\eea
One can verify that $J = X_\xi$ is obtained from \eq{char_sol} with
\[
A(\xi) = \frac{1}{\eta(\xi)},
\quad
B(\tau) =  -\left[\frac{E_+(\tau) - E_-(\tau)}{\eta_+ E_+(\tau) - \eta_- E_-(\tau)}\right].
\]

We then let
\beaa
&&X(\xi,\tau) = \int_0^\xi J\left(\xi^\prime,\tau\right)\, d\xi^\prime,
\\
&&U(\xi,\tau) = X_\tau(\xi,\tau),\quad V(\xi) = F(\xi).
\eeaa
Since $E_-, E_+ > 0$, $\eta_- < \eta_+$, and $\eta_-\le \eta \le \eta_+$, we see from
\eq{defJ} that $J > 0$. It follows that the transformation
$x = X(\xi,\tau)$, $t=\tau$ between spatial and characteristic
coordinates is smoothly invertible, and, according to Proposition~\ref{prop:char_sol},
these expressions define a smooth solution of \eq{uveq1}--\eq{uveq2},
as may be verified directly.
We show that this solution satisfies the required initial and boundary conditions.

First, at $\tau =0$, we have $E_+ = E_- =1$ and $\sigma_+ - \sigma_- = \eta_+ - \eta_-$. It follows from
\eq{defJ} that $J = 1$ at $\tau = 0$, so $x=\xi$, and $v(x,0) = F(x)$.

Second, using the equation
\[
u_x = \frac{U_{\xi}}{X_\xi} = \frac{J_\tau}{J},
\]
we compute from \eq{defJ} that
\[
u_x =  2 \frac{\left(\eta_+-\eta\right)\frac{dE_+}{d\tau}
+ \left(\eta-\eta_-\right)\frac{dE_-}{d\tau}}{\left(\eta_+ - \eta\right)E_+ + \left(\eta - \eta_-\right)
E_-}
- \frac{\frac{d}{d\tau}\left(\sigma_+ - \sigma_-\right)}{\sigma_+ - \sigma_-}.
\]
From \eq{defE+}--\eq{defE-}, we have
\[
\frac{dE_+}{d\tau} = -\frac{1}{4}\left(\sigma_+ - \sigma_-\right) E_+,
\quad
\frac{dE_-}{d\tau} = \frac{1}{4}\left(\sigma_+ - \sigma_-\right) E_-,
\]
and from the jump condition \eq{jumpcons}, we have
\[
\frac{d}{d\tau}\left(\sigma_+ - \sigma_-\right) = -\frac{1}{2}\left(\sigma_+^2 - \sigma_-^2\right).
\]
Using these equations to eliminate $\tau$-derivatives from the expression for $u_x$ and simplifying
the result, we get
\[
u_x = \frac{\left(\eta_+ - \eta\right)\sigma_- E_+ + \left(\eta - \eta_-\right)
\sigma_+ E_-}{\left(\eta_+ - \eta\right)E_+ + \left(\eta - \eta_-\right)
E_-}.
\]
It follows that $u_x=\sigma_-$ at $x=-\infty$, when $\eta = \eta_-$, and
$u_x = \sigma_+$ at $x=\infty$, when $\eta = \eta_+$. \qquad$\Box$

\medskip
For example, let us consider what happens when the derivative $u_x$ vanishes at $x=-\infty$
or $x=\infty$.
If $\sigma_- = 0$, then $\sigma_+ > 0$ satisfies the equation
\[
\frac{d\sigma_+}{dt} + \frac{1}{2} \sigma_+^2 = 0,
\]
which has a global smooth solution forward in time,
\[
\sigma_+(t) = \frac{\sigma_+(0)}{1 + \sigma_+(0) t/2}.
\]
It follows that \eq{uveq1}--\eq{uxcond} has a global smooth solution
forward in time, which may be specified uniquely by the requirement
that $u=v=0$ at $x=-\infty$. This case corresponds to the one that arises
for the fast twist waves analyzed in Section~\ref{sec:fast_twist}.

On the other hand, if $\sigma_+ = 0$, then $\sigma_- < 0$ satisfies the equation
\[
\frac{d\sigma_-}{dt} + \frac{1}{2} \sigma_-^2 = 0,
\]
whose solution
\[
\sigma_-(t) = \frac{\sigma_-(0)}{1 + \sigma_-(0) t/2}
\]
blows up as $t\uparrow t_\ast$, where
$t_\ast = -{2}/{\sigma_-(0)} > 0$.
Thus, a smooth solution of \eq{uveq1}--\eq{uxcond} exists only in
the finite time-interval $0 \le t < t_\ast$.
The derivative $u_x$ blows up simultaneously in the
entire semi-infinite spatial interval to the left of the support
of $v_x$, so it does not appear possible to continue the smooth
solution by any kind of distributional solution after the singularity
forms.

\section{Integrability and Hamiltonian structure}
\label{sec:ham}
\setcounter{equation}{0}

In this section, we show that the twist-wave equation \eq{uveq} is a completely integrable,
bi-Hamiltonian PDE, and that if $(u,v)$ satisfies \eq{uveq} then $u$ satisfies
the HS-equation \eq{HSeq1}.

We begin by describing the relation between \eq{uveq}
and the HS-equation \eq{HSeq1}.
In order to describe
the corresponding relation for the Camassa-Holm (CH) equation at the same time,
we consider the following generalization of \eq{uveq}:
\bea
&&\left(v_t + u v_x\right)_x = 0,
\label{sys1}
\\
&&M u = v_x^2,
\label{sys2}
\eea
where $M$ is a self-adjoint linear operator
acting on functions of $x$ that commutes with $\partial_x$.
If $M = \partial_x^2$, then \eq{sys1}--\eq{sys2} is \eq{uveq}.

We suppose that $u$, $v$ are smooth solutions of \eq{sys1}--\eq{sys2}.
Differentiating \eq{sys2} with respect to $t$, using \eq{sys1} to write $v_{xt}$
in terms of $u$, $v$ and their spatial derivatives, then
using \eq{sys2} to eliminate $v$ from the result, we find that $u$ satisfies
\be
m_t + m u_x + \left(m u\right)_x = 0\qquad \mbox{with $m = Mu$}.
\label{genHS}
\ee
If $M=\partial_x^2$, then \eq{genHS} is the HS-equation \eq{HSeq1};
if $M = \partial_x^2 - 1$, then \eq{genHS} is the CH-equation \cite{CH},
\[
\left[\left(u_t+uu_x\right)_x -\frac{1}{2} u_x^2\right]_x = u_t +3 u u_x.
\]
Conversely, if $m > 0$ is a smooth solution of \eq{genHS} and
$v_x = \sqrt{m}$, then $v$ satisfies \eq{sys1}--\eq{sys2}.

Because of the nonlinearity of this transformation, difficulties may arise in
its application to distributional solutions. For example, the function
\[
u(x,t) =  \left\{\begin{array}{ll} 0 & \mbox{if $x\le 0$}
\\
                                  2x/t& \mbox{if $x>0$}
                            \end{array}\right.
\]
is a weak solution of the HS-equation \eq{HSeq1} in $t > 0$, and
\[
u_{xx}(x,t) = \frac{2}{t}\delta(x)
\]
is non-negative in the sense of distributions.
There is, however, no standard way to define a distribution
$v$ such that $u_{xx}=v_x^2$.

The HS-equation \eq{HSeq1} is bi-Hamiltonian and completely integrable \cite{hz},
so \eq{uveq} is also.
Next, we consider the effect of the transformation $v\mapsto u$ on the
Hamiltonian structures of these equations.

The system \eq{sys1}--\eq{sys2} is obtained from the variational principle
\[
\delta \int \frac{1}{2}\left\{-v_t v_x - u v_x^2 + \frac{1}{2}u M u\right\}\, dx dt = 0.
\]
Variations with respect to $u$ yield \eq{sys2}, and variations with respect to $v$
yield \eq{sys1}.
We may eliminate $u$ by means of the constraint
equation $u = M^{-1} \left(v_x^2\right)$ to obtain a variational principle for $v$ alone,
\be
\delta \int \frac{1}{2}\left\{- v_t v_x - \frac{1}{2} v_x^2M^{-1} \left(v_x^2\right)
\right\}\, dx dt = 0.
\label{action2}
\ee
The Euler-Lagrange equation for \eq{action2},
\be
\left[v_t + M^{-1} \left(v_x^2\right)v_x\right]_x = 0,
\label{onlyveq}
\ee
is equivalent to \eq{sys1}--\eq{sys2}.
Here, and below, we assume that operators such as $M$ and $\partial_x$
are invertible; in the case of differential operators,
this requires the addition of suitable boundary conditions
which we do not specify explicitly.

Making a Legendre transform of the Lagrangian in \eq{action2}, we get the
corresponding Hamiltonian form of \eq{onlyveq},
\be
v_t = \partial_x^{-1} \left(\frac{\delta\mathcal{H}}{\delta v}\right),
\qquad \mathcal{H} = \frac{1}{4} \int v_x^2M^{-1} \left(v_x^2\right)\, dx
\label{hamveq},
\ee
where $\partial_x^{-1}$ is
the constant Hamiltonian operator associated canonically with the
variational principle \eq{action2}.

\begin{proposition}
\label{prop:hamtrans}
Let $\{\mathcal{F},\mathcal{G}\}$ denote the Poisson bracket of functionals
$\mathcal{F}$, $\mathcal{G}$ of $v$ associated with
the constant Hamiltonian operator $\partial_x^{-1}$,
\be
\left\{\mathcal{F},\mathcal{G}\right\} =
\int \frac{\delta \mathcal{F}}{\delta v}\, \partial_x^{-1}
\left(\frac{\delta \mathcal{G}}{\delta v}\right)\, dx.
\label{canpb}
\ee
Under the change of variables $u = M^{-1}\left(v_x^2\right)$,
where $M$ is a self-adjoint linear operator,
the bracket \eq{canpb} transforms formally
into a Lie-Poisson bracket
\bea
&&\left\{\mathcal{F},\mathcal{G}\right\} =
\int \frac{\delta \mathcal{F}}{\delta u}\,J(u)
\left(\frac{\delta \mathcal{G}}{\delta u}\right)\, dx,
\label{lppb1}
\\
&&J = -2 M^{-1}\left(m\partial_x + \partial_x m\right)M^{-1},
\label{lppb2}
\eea
where $m = Mu$.
\end{proposition}

\medskip\noindent
\textbf{Proof.} First, we consider the nonlinear change of
variables $m=v_x^2$. Variations of $v$ of the form
$v^\eps = v + \eps k + \dots$ lead to variations
$m^\eps = m + \eps h + \dots$  of $m$ where
$h = 2 v_x k_x$.
For any functional $\mathcal{F}$ of $m$, with
$\mathcal{F}^\eps = \mathcal{F}(m^\eps)$, we have
\[
\left.\frac{d}{d\eps}\mathcal{F}^\eps\right|_{\eps=0} =
\int \frac{\delta \mathcal{F}}{\delta m} \,h\,dx.
\]
Writing $h$ in terms of $k$ and using the skew-adjointness of $\partial_x$, we compute that
\[
\left.\frac{d}{d\eps}\mathcal{F}^\eps\right|_{\eps=0}
= - 2\int \left(\frac{\delta \mathcal{F}}{\delta m} v_x\right)_x k \,dx.
\]
Since
\[
\left.\frac{d}{d\eps}\mathcal{F}^\eps\right|_{\eps=0} =
\int \frac{\delta \mathcal{F}}{\delta v} \,k\,dx,
\]
we conclude that
\[
\frac{\delta \mathcal{F}}{\delta v} = - 2 \left(\frac{\delta \mathcal{F}}{\delta m}
 \sqrt{m}\right)_x.
\]
Using this equation in \eq{canpb} and integrating by parts, we get
\beaa
\left\{\mathcal{F},\mathcal{G}\right\} &=&
4\int \left(\frac{\delta \mathcal{F}}{\delta m}
 \sqrt{m}\right)_x \frac{\delta \mathcal{G}}{\delta m}
 \sqrt{m}\, dx
\\
&=& - 2\int \frac{\delta \mathcal{F}}{\delta m}
 \left(m\partial_x + \partial_x m\right)\frac{\delta \mathcal{G}}{\delta m}
\, dx.
\eeaa
Making the linear change of variables $u = M^{-1} m$ in this expression, and using the self-adjointness
of $M$, we get \eq{lppb1}--\eq{lppb2}.
\qquad$\Box$

\medskip
The Hamiltonian form of \eq{genHS} for $u$ corresponding to \eq{hamveq} for $v$ is therefore
\[
u_t = J\left( \frac{\partial \mathcal{H}}{\partial u}\right),
\qquad
\mathcal{H} = \frac{1}{4} \int u Mu \, dx,
\]
where $J$ is given in \eq{lppb2}.

To give a second Hamiltonian structure for \eq{sys1}--\eq{sys2},
we define a skew-adjoint operator
$K$, depending on $v$, by
\be
K = v_x \partial_x^{-1} M^{-1} v_x,
\label{defK}
\ee
where $M$ is a self-adjoint linear
operator commuting with $\partial_x$,
as before.

We find that the operator \eq{defK} satisfies the Jacobi identity if the quantity
\[
\left(g h_x - h g_x \right)M f_x +  \left(h f_x - f h_x \right)M g_x
+ \left(f g_x - g f_x \right)M h_x
\]
is an exact $x$-derivative for arbitrary functions $f$, $g$, $h$.
This condition holds for $M = \partial_x^2$, since
\[
 \left(g h_x - h g_x \right)f_{xxx} =\left[\left(g h_x-hg_x\right)f_{xx}\right]_x
+ h f_{xx} g_{xx} - g h_{xx} f_{xx},
\]
and the terms that are not exact derivatives cancel under a cyclic summation.
The condition also holds for $M=\partial_x^2 -1$. Moreover, in those cases,
$\left(c_1 \partial_x^{-1} + c_2 K\right)$ satisfies the Jacobi identity for arbitrary
real constants $c_1$, $c_2$, so that $\partial_x^{-1}$ and $K$ define
compatible Hamiltonian structures.

The Hamiltonian form of \eq{onlyveq} with respect to $K$ is
\be
v_t = K\left(\frac{\delta\mathcal{P}}{\delta v}\right),
\qquad\mathcal{P} = \int v_x^2 \, dx.
\label{hamveq1}
\ee
If $M=\partial_x^2$, then $K = v_x \partial_x^{-3} v_x$,
and \eq{hamveq1} is equivalent to \eq{uveq}.

Under the transformation $m=v_x^2$, equation \eq{hamveq1} becomes
\beaa
&&m_t = \widetilde K \left(\frac{\delta\mathcal{P}}{\delta m}\right),
\qquad
\mathcal{P} = \int m \, dx,
\\
&&\widetilde{K} = - \left(m\partial_x + \partial_x m\right) \left(\partial_x^{-1} M^{-1}\right)
\left(m\partial_x + \partial_x m\right),
\eeaa
which gives \eq{genHS}.

When $M = \partial_x^2$, equations \eq{hamveq} and \eq{hamveq1} provide
a bi-Hamiltonian structure for \eq{uveq}. One can then obtain
an infinite sequence of commuting Hamiltonian
flows by recursion. We will not write them out explicitly here,
but we remark that among them is a Hamiltonian structure
for $v$ which maps to the Hamiltonian structure for $u$ canonically
associated with the variational principle in \eq{hsaction}.

Lax pairs for \eq{uveq} follow directly by transformation of the Lax pairs
for the HS-equation \cite{hz}. We define
\[
L = \partial_x^{-1} v_x^2\partial_x^{-1},
\qquad
A = \frac{1}{2} \left(u\partial_x + \partial_x u\right).
\]
Then, using the identity
$f\partial_x^{-1} - \partial_x^{-1} f = \partial_x^{-1} f_x \partial_x^{-1}$,
we compute that the Lax equation
$L_t = \left[L,A\right]$ is equivalent to
\[
\left(v_x^2\right)_t + \left(u v_x^2 + \frac{1}{2} u_x^2\right)_x = 0,
\qquad
u_{xx} = v_x^2,
\]
which may be rewritten as \eq{uveq}.
Alternatively, we can set $u = \partial_x^{-2} \left(v_x^2\right)$
in the original Lax pair.



\section{Periodic twist waves}
\label{sec:per}
\setcounter{equation}{0}

Spatially periodic splay waves are described by the following version of the HS-equation
\cite{per,rev}
\[
\left(u_{t} +  u u_x\right)_x = \frac{1}{2}\left(u_x^2 - \langle u_x^2\rangle\right).
\]
Here, $u(x,t)$ is a periodic function of $x$, and
angular brackets denote an average over a period.
The wave also drives a mean-field, which evolves on the same
time-scale, $t=O(1)$, as the wave.

The interaction between a weakly nonlinear twist wave and a mean-field
is more complicated  because the mean-field evolves on a
faster time-scale than the wave. This is a consequence of the fact
that the nonlinear self-interaction of
a weakly nonlinear twist wave is cubic, but the mean-field is driven
by quadratic nonlinearities.

In this section, we derive the following generalization of the twist-wave
asymptotic equation \eq{uveq} that applies to periodic waves:
\bea
&&\left(v_{t} +  u v_x\right)_x + \mu \langle v_x^2\rangle v = 0,
\label{uvpereq1}
\\
&&
u_{xx} = v_x^2 - \langle v_x^2\rangle.
\label{uvpereq2}
\eea
Here, $u(x,t)$, $v(x,t)$ are periodic functions of $x$, which
we assume to have zero mean without loss generality,
and $\mu$ is a constant that cannot be removed by rescaling.
We will derive \eq{uvpereq1}--\eq{uvpereq2} from the one-dimensional
wave equations \eq{wavea}--\eq{waveb},
but a similar derivation would apply to more general systems.

It is interesting to note that the mean-field interaction introduces a
dispersive term of Klein-Gordon type into the evolution equation
\eq{uvpereq1} for $v$,
despite the fact that the original system is scale-invariant and non-dispersive.
The scale-invariance is preserved by the fact that the coefficient
of the dispersive term is proportional to the mean wave-energy (or momentum)
$\langle v_x^2\rangle$, which is a constant conserved quantity.

The mean-terms prevent the elimination of $v$ from the system \eq{uvpereq1}--\eq{uvpereq2}
by cross-differentiation, as is possible in the case of \eq{sys1}--\eq{sys2}.
Instead, one finds that
\[
\left[\left(u_t + u u_x\right)_x - \frac{1}{2}u _x^2
- \langle v_x^2\rangle \left(2 u + \mu  v^2\right)\right]_x = 0,
\]
which suggests that \eq{uvpereq1}--\eq{uvpereq2} is not completely integrable
when $\mu\ne 0$.

\subsection{Derivation of the periodic equation}

We consider spatially-periodic solutions of \eq{wavea}--\eq{waveb}
with period of the order $\eps$ in which
$\psi$ has amplitude of the order $\eps^{1/2}$,
where $\eps$ is a small parameter. The corresponding time-scale
for the nonlinear evolution of the $\psi$-wave is of the order $1$. As we will see,
the $\psi$-wave generates a mean $\phi$-field which
evolves on a time-scale of the order $\eps^{1/2}$. This
mean field modulates the speed of
the $\psi$-wave on the same time-scale.

We therefore introduce multiple-scale variables
\bea
&&\theta=\frac{x-\eps^{1/2}s\left({t}/{\eps^{1/2}}\right)}{\eps},
\nonumber
\\
&&\tau = \frac{t}{\eps^{1/2}},
\label{msper}
\\
&&t = t,
\nonumber
\eea
where $s(\tau)$ is a suitable phase function, and
look for an asymptotic solution of \eq{wavea}--\eq{waveb} of the form
\beaa
&&\phi = \phi_0\left(\tau\right)
+ \eps \phi_2\left(\theta,\tau,t\right) + O(\eps^{3/2}),
\\
&&\psi = \eps^{1/2}\psi_1\left(\theta,\tau,t\right)
+ \eps \psi_2\left(\theta,\tau,t\right) + O(\eps^{3/2}),
\eeaa
where all the terms are periodic functions of $\theta$. We will also require
below that the terms are periodic functions of $\tau$.

We use this ansatz in \eq{wavea}--\eq{waveb}, expand derivatives as
\[
\partial_t \to - \frac{1}{\eps} s_\tau\partial_\theta + \frac{1}{\eps^{1/2}} \partial_\tau + \partial_t,
\quad \partial_x \to \frac{1}{\eps} \partial _\theta,
\]
Taylor expand the result with respect to $\eps$, and equate coefficients of powers
of $\eps^{1/2}$ to zero.

We consider first the $\psi$-equation \eq{waveb}.
At the order $\eps^{-3/2}$, we obtain that
\[
\left(s_\tau^2 - b_0^2\right) \psi_{1\theta\theta} = 0,
\]
where the $0$-subscript on a function of $\phi$ denotes evaluation at $\phi=\phi_0$.
It follows from this equation that $s_\tau^2 = b_0^2$.
For definiteness, we consider a right-moving wave and assume that
\be
s_\tau = b_0,
\label{stau}
\ee
where $b_0 > 0$.

At the order $\eps^{-1}$, we obtain that
\be
2 s_\tau \psi_{1\theta\tau} + \left(s_{\tau\tau} + \frac{2q_{0\tau}}{q_0} s_\tau\right) \psi_{1\theta} = 0.
\label{psi11eq}
\ee
We may choose $\psi_1$ so that it is a zero-mean periodic function of $\theta$.
It follows from \eq{stau} and \eq{psi11eq} that
\be
\psi_1(\theta,\tau,t) = \frac{v(\theta,t)}{q_0(\tau) b_0^{1/2}(\tau)}
\label{psi1eq}
\ee
where $v(\theta,t)$ is a zero-mean periodic function of $\theta$ which is independent of $\tau$.

At the order $\eps^{-1/2}$, we obtain that
\[
2 s_\tau \psi_{2\theta\tau} + \left(s_{\tau\tau} + \frac{2q_{0\tau}}{q_0} s_\tau\right) \psi_{2\theta}
+ 2 s_\tau \psi_{1\theta t} + \left(2 b_0 b_0^\prime \phi_2 \psi_{1\theta}\right)_\theta
- \psi_{1\tau\tau} = 0.
\]
Using \eq{stau}, we may write this equation as
\be
\left(q_0 b_0^{1/2} \psi_{2\theta} \right)_\tau
 + q_0 b_0^{1/2} \psi_{1\theta t}
+ \left(q_0 b_0^{1/2} b_0^\prime\phi_2 \psi_{1\theta}\right)_\theta
-\frac{q_0}{2 b_0^{1/2}} \psi_{1\tau\tau} = 0.
\label{psi2eq}
\ee
We will return to \eq{psi2eq} after we expand the $\phi$-equation.

The leading-order terms in the expansion of the $\phi$-equation \eq{wavea}
are of the order $\eps^{-1}$, and give
\[
\left(s_\tau^2 - a_0^2\right)\phi_{2\theta\theta} + q_0^2 b_0 b_0^\prime \psi_{1\theta}^2 + \phi_{0\tau\tau}
= 0.
\]
Using \eq{stau} and \eq{psi1eq} in this equation, we get
\be
\left(b_0^2 - a_0^2\right)\phi_{2\theta\theta} + b_0^\prime v_{\theta}^2 + \phi_{0\tau\tau}
\label{phizeqprime}
= 0.
\ee
Averaging this equation with respect to $\theta$, we find that
\be
\phi_{0\tau\tau} +\langle v_\theta^2\rangle  b_0^\prime = 0,
\label{phizeq}
\ee
where the angular brackets denote an average with respect to $\theta$.
As we will see, for smooth solutions, the quantity $\langle v_\theta^2\rangle$
is a constant independent of $t$, so \eq{phizeq}
is consistent with the ansatz that $\phi_0$ depends only on $\tau$.

Equation \eq{phizeq} provides an ODE for $\phi_0$, corresponding to motion in a potential
proportional to the twist-wave speed $b_0=b(\phi_0)$. For definiteness,
we assume that the solution of \eq{phizeq} for $\phi_0(\tau)$ is a periodic
function of $\tau$. We then require
that all other terms in the expansion are periodic functions of $\tau$.

Subtracting \eq{phizeq} from \eq{phizeqprime},
we find that
\[
\left(b_0^2 - a_0^2\right) \phi_{2\theta\theta}
+ b_0^\prime \left(v_{\theta}^2 - \langle v_{\theta}^2\rangle\right) = 0.
\]
Hence, we may write
\be
\phi_2(\theta,\tau,t) = \left[\frac{b_0^\prime(\tau)}{a_0^2(\tau) - b_0^2(\tau)}\right] u(\theta,t),
\label{phi2eq}
\ee
where $u$ satisfies
\[
u_{\theta\theta} = v_{\theta}^2 - \langle v_{\theta}^2\rangle.
\]

Using \eq{psi1eq} and \eq{phi2eq} in \eq{psi2eq}, we get
\[
\left(q_0 b_0^{1/2} \psi_{2\theta} \right)_\tau
 + v_{\theta t}
+ \left[\frac{\left(b_0^\prime\right)^2}{a_0^2 - b_0^2}\right]\left(u v_{\theta}\right)_\theta
-\frac{q_0}{2 b_0^{1/2}} \left(\frac{1}{q_0 b_0^{1/2}}\right)_{\tau \tau} v = 0.
\]
Averaging this equation with respect to $\tau$, we get
\[
\left(v_{t} + \Lambda u v_\theta\right)_\theta + N v = 0,
\]
where
\beaa
\Lambda &=& \oint \frac{\left(b_0^\prime\right)^2}{a_0^2 - b_0^2} \, d\tau,
\\
N &=& - \frac{1}{2}\oint \frac{q_0}{b_0^{1/2}} \left(\frac{1}{q_0 b_0^{1/2}}\right)_{\tau \tau}\, d\tau.
\eeaa
Here,  $\oint$ denotes an average over a period in $\tau$.

To make the dependence of $\phi_0$ on $v$ explicit, we introduce a new time variable
\[
T = \langle v_\theta^2\rangle^{1/2} \tau.
\]
We may then rewrite \eq{phizeq} as
\be
\phi_{0TT} +  b_0^\prime = 0.
\label{phizode}
\ee
Given a solution of this equation for $\phi_0(T)$, we have
$N = \langle v_\theta^2\rangle M$,
where
\[
M = - \frac{1}{2}\oint \frac{q_0}{b_0^{1/2}} \left(\frac{1}{q_0 b_0^{1/2}}\right)_{TT}\, d T
\]
is a constant independent of $v$, and $\oint$ denotes an average with respect to $T$ over a period.
From \eq{phizode}, we have
\[
\frac{1}{2}\phi_{0T}^2 + b_0 = E
\]
for some constant $E$. Using an integration by parts, we may also write $M$ as
\beaa
M &=&
\frac{1}{2} \oint \frac{1}{b_0} \left(\frac{b_{0T}^2}{4 b_0^2} - \frac{q_{0T}^2}{q_0^2}\right)\, dT
\\&=& \oint \frac{E-b_0}{b_0} \left[\frac{\left(b_{0}^\prime\right)^2}{4 b_0^2} - \frac{\left(q_{0}^\prime\right)^2}{q_0^2}\right]\, dT.
\eeaa

Thus, the final equations for $u(\theta,t)$, $v(\theta,t)$ are
\bea
&&\left(v_{t} + \Lambda u v_\theta\right)_\theta + M \langle v_\theta^2\rangle v = 0,
\label{asyper1}
\\
&&u_{\theta\theta} = v_\theta^2 - \langle v_\theta^2\rangle.
\label{asyper2}
\eea
It follows from these equations that, for smooth solutions,
\[
\left(v_\theta^2\right)_t + \left[\Lambda \left(u v_\theta^2 - \frac{1}{2} u_\theta^2
+ \langle v_\theta^2\rangle u_\theta\right)
+ M \langle v_\theta^2\rangle v^2\right]_\theta = 0.
\]
Taking the average of this equation with respect to $\theta$,
we find that $\langle v_\theta^2\rangle$ is constant in time,
as stated earlier.

In summary, the asymptotic solution of \eq{wavea}--\eq{waveb} is given by
\beaa
\phi &=& \phi_0(\tau) + \frac{\eps b_0^\prime(\tau)}{a_0^2(\tau)-b_0^2(\tau)}\, u(\theta,t) +O(\eps^{3/2}),
\\
\psi &=&\frac{\eps^{1/2}}{q_0(\tau) b_0^{1/2}(\tau)}\, v(\theta,t) + O(\eps),
\eeaa
where the multiple-scale variables $\theta$, $\tau$ are evaluated at \eq{msper}, $s$ satisfies \eq{stau},
$\phi_0$ satisfies \eq{phizeq}, and $u$, $v$ satisfy \eq{asyper1}--\eq{asyper2}.

If $\Lambda \ne 0$ then we may rescale variables in \eq{asyper1}--\eq{asyper2}
to get \eq{uvpereq1}--\eq{uvpereq2} with
\[
\mu = \frac{M}{\Lambda}.
\]

As an example, let us consider the wave speeds in \eq{dircfs} arising from the
one-dimensional director field equations, where
\[
b(\phi) = \sqrt{\beta \sin^2 \phi + \gamma \cos^2\phi}.
\]
If $\beta < \gamma$, then $b$ has a minimum at $\phi=\pi/2$, and \eq{phizeq}
has periodic solutions for the mean field $\phi_0$ that oscillate around $\pi/2$.
Our asymptotic solution applies in this case.
If $\beta > \gamma$, then $b$ has a minimum at $\phi=0$. Although \eq{phizeq}
also has periodic solutions in this case, there is a
loss of strict hyperbolicity at $\phi=0$, where
$a=b$, and the asymptotic solution breaks down.


\bigskip\noindent
\textbf{Acknowledgements.} The work of J.K.H. was partially
supported by the NSF under grant number DMS--0607355.


\appendix

\section{Algebraic details}
\setcounter{equation}{0}
\label{app:details}

\subsection{Linearized equations}
\label{app:lin}

The linearization of \eq{dirPDE} at $\bnz$ is
\be
\bnp_{tt} = \alpha\nabla(\div \bnp) - \beta \curl\left(\ap\bnz\right)
- \gamma\curl\left(\bBp\times\bnz\right) + \lambdap\bnz,
\label{lin}
\ee
where
\[
\ap = \bnz\cdot\curl\bnp,\qquad
\bBp = \bnz\times \curl\bnp,
\]
and the Lagrange multiplier $\lambdap$ is chosen so that
$\bnz\cdot \bnp = 0$.

The Fourier mode
\[
\bnp(\bx,t) = \bnh \, e^{i\bk\cdot\bx - i\omega t},
\qquad
\lambdap(\bx,t) = \lambdah \, e^{i\bk\cdot\bx - i\omega t}
\]
is a solution of the linearized equations \eq{lin} if
\be
\mathcal{L}\bnh - \lambdah \bnz = 0,
\qquad \bnz\cdot\bnh = 0,
\label{linevp}
\ee
where the linear map $\mathcal{L} : \Rl^3 \to \Rl^3$ is defined by
\bea
&&\mathcal{L}\bnh = \omega^2\bnh - \alpha\left(\bk\cdot\bnh\right)\bk
+\beta \tA \left(\bk\times\bnz\right)
+ \gamma\bk\times\left(\tB\times\bnz\right),\quad
\label{defL}
\\
&&\qquad\tA = \bnz\cdot\left(\bk\times\bnh\right),\qquad
\tB = \bnz\times\left(\bk\times\bnh\right).
\nonumber
\eea
We solve this eigenvalue problem in the next proposition.

\begin{proposition}
\label{prop:linsys}
Suppose that $\bnz, \bk\in \Rl^3$ where
$\bnz$ is a unit vector and $\bk$ is non-zero,
and $\alpha, \beta, \gamma \in \Rl$ are distinct positive constants.
For $\omega\in \Rl$
let $\mathcal{L} : \Rl^3 \to \Rl^3$ be the linear map defined in
\eq{defL}. Then the linear system
\be
\mathcal{L} \bm - \lambda \bn_0 = \bF,
\qquad\bn_0\cdot \bm = G
\label{lineq}
\ee
has a unique solution for $\{\bm,\lambda\}\in \Rl^3\times \Rl$
for every $\{\bF,G\}\in \Rl^3\times \Rl$ unless
$\omega^2=\Wsplay^2(\bk;\bnz)$ or $\omega^2=\Wtwist^2(\bk;\bnz)$,
where $\Wsplay^2$, $\Wtwist^2$ are defined in \eq{amode}--\eq{bmode}.

\smallskip\noindent
(a) If $\bk$ is not parallel to $\bnz$ and $\omega^2 = \Wsplay^2(\bk;\bnz)$, then
the general solution of \eq{lineq} when $\bF=0$, $G=0$ is
\be
\bm = c\Rsplay,\qquad
\lambda = - c\left(\alpha -\gamma\right)\left(\bk\cdot\bnz\right)\left(\bk\cdot\Rsplay\right),
\label{homsplaysol}
\ee
where $c$ is an arbitrary constant and  $\Rsplay = \bk - \left(\bk\cdot\bnz\right)\bnz$.
Equation \eq{lineq} is solvable for $\{\bm, \lambda\}$
if and only if $\{\bF,G\}$ satisfy
\be
\Rsplay\cdot \bF + \left(\alpha-\gamma\right)\left(\bk\cdot\bnz\right)\left(\bk\cdot\Rsplay\right)G = 0.
\label{solsplay}
\ee
(b) If $\bk$ is not parallel to $\bnz$ and $\omega^2 = \Wtwist^2(\bk;\bnz)$, then the general solution
of \eq{lineq} when $\bF=0$, $G=0$ is
\be
\bm = c\Rtwist,\qquad \lambda = 0,
\label{homtwistsol}
\ee
where $c$ is an arbitrary constant and $\Rtwist = \bk \times \bnz$.
Equation \eq{lineq}
is solvable for $\{\bm, \lambda\}$ if and only if $\bF$ satisfies
\be
\Rtwist\cdot \bF = 0.
\label{soltwist}
\ee
(c) If $\bk$ is parallel to $\bnz$ and $\omega^2 = \gamma(\bk\cdot\bnz)^2$, then
the general solution of \eq{lineq} when $\bF=0$, $G=0$ is $\{\bm,\lambda\} = \{\bm^\perp,0\}$
where $\bm^\perp$ is any vector
orthogonal to $\bnz$. Equation
\eq{lineq} is solvable for $\{\bm, \lambda\}$ if and only if $\bF$ is parallel
to $\bnz$.
\end{proposition}

\medskip\noindent
\textbf{Proof.} First, we suppose that $\bk$ is not parallel to $\bnz$.
Expanding
\beaa
&&\bm = m_1 \bk + m_2 \Rsplay + m_3 \Rtwist,
\\
&&\bF = F_1 \bk + F_2 \bnz + F_3 \Rtwist,
\eeaa
we find, after some algebra, that \eq{lineq} is equivalent to
\beaa
&&\left(\bk\cdot\bnz\right) m_1 = G,
\\
&&\lambda + \left(\bk\cdot\bnz\right) \left(\omega^2 - \gamma k^2\right) m_2 =  -F_2,
\\
&&\left(\omega^2 - \alpha k^2\right) m_1 + \left(\omega^2 - \Wsplay^2\right) m_2 = F_1,
\\
&&\left(\omega^2-\Wtwist^2\right) m_3 = F_3.
\eeaa
The first two equations determine $m_1$ and $\lambda$, and the remaining two equations
determine $m_2$, $m_3$ unless $\omega^2=\Wsplay^2$ or $\omega^2=\Wtwist^2$.

If $\omega^2=\Wsplay^2$, then $\omega^2\ne b^2$
(since $\bk$ is not parallel to $\bnz$) and
the fourth equation is solvable
for $m_3$. The third equation is solvable for $m_2$ if and only if
\[
 \left(\bk\cdot\bnz\right)F_1 - \left(\Wsplay^2 - \alpha k^2\right)G = 0.
\]
Using the equations
\[
\Wsplay^2-\alpha k^2 = -\left(\alpha -\gamma\right)\left(\bk\cdot\bnz\right)^2,\quad
\Rsplay\cdot \bF = \left(\bk\cdot\Rsplay\right) F_1
\]
in this condition, we get \eq{solsplay}. If $\bF=0$, $G=0$, then we find that $m_1=0$, $m_3=0$
and $m_2=c$, where $c$ is an arbitrary constant. Computing the corresponding value
of $\lambda$, we get \eq{homsplaysol}.

If $\omega^2 = \Wtwist^2$, then the first three equations are solvable
for $\lambda$, $m_1$, $m_2$. The last equation is solvable for $m_3$ if and only if
$F_3=0$, which gives \eq{soltwist}. If $\bF=0$, $G=0$ then $\lambda=0$, $m_1=0$,
$m_2=0$, and $m_3=c$, which gives \eq{homtwistsol}.

Finally, we suppose that $\bk = k\bnz$ is parallel to $\bnz$. Then
\[
\mathcal{L}\bm  = \left(\omega^2 -\gamma k^2\right)\bm
-(\alpha-\gamma)k^2 \left(\bnz\cdot \bm\right) \bnz.
\]
Equation \eq{lineq}
is therefore uniquely solvable unless $\omega^2 = \gamma k^2$,
when
\[
\mathcal{L}\bm  = -(\alpha-\gamma)k^2 \left(\bnz\cdot \bm\right) \bnz.
\]
In that case, \eq{lineq}
is solvable if and only if $\bF = F \bnz$ is parallel to $\bnz$, and the solution is

\[
\bm = G \bnz + \bm^\perp,\qquad \lambda = -\left[F + (\alpha-\gamma)k^2 G\right],
\]
where $\bm^\perp$ is an arbitrary vector orthogonal to $\bnz$.
$\Box$

\medskip
From Proposition~\ref{prop:linsys}, the
solutions of the eigenvalue problem \eq{linevp} are given by
\eq{lindisp}--\eq{Rtwist}.

\subsection{Weakly nonlinear splay waves}
\label{app:splay}

We look for an asymptotic solution of \eq{dirPDE} of the form
\eq{nsplayexp}--\eq{nsplayexp1}. We expand derivatives as
\be
\partial_t \to -\frac{\omega}{\eps}\partial_\theta + \partial_t,
\quad
\nabla \to \frac{\bk}{\eps}\partial_\theta + \nabla,
\label{derivexp}
\ee
where $\omega$, $\bk$ are defined in \eq{locfreq}.
The corresponding expansions of $\lambda$, $A$, $\bB$ are
\[
\lambda = \frac{1}{\eps}\lambda_1 + \lambda_2  + \dots,
\quad
A = A_1 + \eps A_2  + \dots,
\quad
\bB = \bB_1 + \eps\bB_2  + \dots.
\]
where
\beaa
&&A_1 = \bnz\cdot\left(\bk\times\bn_{1\theta}\right)+ \bnz\cdot \curl\bnz ,
\\
&&\bB_1 = \bnz\times\left(\bk\times\bn_{1\theta}\right) +\bnz\times \curl\bnz,
\\
&&A_2 = \bnz\cdot\left(\bk\times\bn_{2\theta}\right)+\bn_1\cdot\left(\bk\times\bn_{1\theta}\right)
+ \bnz\cdot\curl \bn_1 + \bn_1\cdot\curl\bnz ,
\\
&&\bB_ 2 = \bnz\times\left(\bk\times\bn_{2\theta}\right)+ \bn_1\times\left(\bk\times\bn_{1\theta}\right)
 + \bnz\times\curl \bn_1 + \bn_1\times\curl\bnz.
\eeaa
We write these expressions as
\bea
&&A_i = \tA_i + \rA_i,
\qquad\qquad\quad
\bB_i =\tB_i + \rB_i,
\label{Ahat}
\\
&&\tA_i = \bnz\cdot\left(\bk\times\bn_{i\theta}\right),
\qquad
\tB_i = \bnz\times\left(\bk\times\bn_{i\theta}\right).
\nonumber
\eea

We use these expansions in \eq{dirPDE}, Taylor expand the result with respect
to $\eps$, and equate coefficients of powers of $\eps$.
At the order $\eps^{-1}$, we obtain
linearized equations for $\bn_1$, $\lambda_1$, which have the form
\be
\mathcal{L} \bn_{1\theta\theta} - \lambda_1 \bnz =0,
\qquad
\bnz\cdot \bn_1 = 0,
\label{linsysasy}
\ee
where $\mathcal{L}$ is the linear map defined
in \eq{defL}.

The system \eq{linsysasy} has the splay-wave eigenvalue
$\omega^2 = \Wsplay^2\left(\bk;\bnz\right)$, where $\Wsplay$ is defined
in \eq{amode}. The corresponding solution for $\bn_{1}$, after two integrations
with respect to $\theta$, is
\[
\bn_1(\theta,\bx,t) = u(\theta,\bx,t) \Rsplay(\bx,t),
\]
where $u$ is an arbitrary scalar-valued function, and $\Rsplay$ is defined in
\eq{Rsplay}. The solution for $\lambda_1$ is
\[
\lambda_1 = - \left(\alpha-\gamma\right)u_{\theta\theta} \left(\bk\cdot\bnz\right)
\left(\bk\cdot\Rsplay\right).
\]

At the order $\eps^0$, we obtain equations for $\bn_2$, $\lambda_2$.
Using the fact that $\bnz$ is a solution of \eq{dirPDE}, with $\lambda=\lambda_0$ say,
we may
write them as
\be
\mathcal{L}\bn_{2\theta\theta} - \rlambda_2 \bn_0 = \bF_1,
\qquad
\bnz\cdot\bn_2 = G_1,
\label{linsysasy1}
\ee
where $\rlambda_2 = \lambda_2 - \lambda_0$ and
\beaa
&&\bF_1 = 2\omega \bn_{1\theta t} + \omega_t \bn_{1\theta} + \alpha\left\{\left(\div \bn_{1\theta}
\right)\bk + \nabla\left(\bk\cdot \bn_{1\theta}\right)\right\}
\\
&&\qquad\quad-\beta\left\{\rA_{2\theta}\left(\bk\times\bnz\right) + A_1\left(\bk\times\bn_{1\theta}\right)
+ \bk\times\left(A_1\bn_1\right)_\theta\right.
\\
&&\qquad\qquad\qquad\qquad \left.+ \tA_1\curl\bnz + \curl\left(\tA_1\bnz\right)\right\}
\\
&&\qquad\quad+ \gamma\left\{- \bk \times \left(\rB_{2\theta}\times\bnz \right)
+ \bB_1\times\left(\bk\times\bn_{1\theta}\right)- \bk\times\left(\bB_1 \times \bn_1\right)_\theta
\right.
\\
&&\qquad\qquad\qquad\qquad \left. + \tB_1\times \curl\bnz
-\curl\left(\tB_1\times \bnz\right) \right\} + \lambda_1 \bn_1,
\\
&&G_1 = - \frac{1}{2} \bn_1\cdot\bn_1.
\eeaa
From Proposition~\ref{prop:linsys},
equation \eq{linsysasy1} is solvable for $\bn_{2\theta\theta}$ and $\rlambda_2$ if and only if
\[
\Rsplay\cdot\bF_1 + \left(\alpha-\gamma\right) \left(\bk\cdot\bnz\right)
\left(\bk\cdot\Rsplay\right)G_{1\theta\theta} = 0.
\]
After some algebra, we find that this solvability condition
gives \eq{hseq}.

\subsection{Weakly nonlinear twist waves}
\label{app:twist}

We look for an asymptotic solution of \eq{dirPDE} of the form
\eq{nexpeps}--\eq{nexp}, expanding derivatives as in \eq{derivexp}.
The corresponding expansions of $\lambda$, $A$, $\bB$ are
\beaa
&&\lambda = \frac{1}{\eps^{3/2}}\lambda_1 + \frac{1}{\eps}\lambda_2
 + \frac{1}{\eps^{1/2}}\lambda_3 + \dots,
\\
&&A = \frac{1}{\eps^{1/2}} A_1 + A_2 +  \eps^{1/2}A_3 + \dots,
\\
&&\bB = \frac{1}{\eps^{1/2}} \bB_1 +  \bB_2 + \eps^{1/2}\bB_3 + \dots,
\eeaa
where
\beaa
&&A_1 = \bnz\cdot\left(\bk\times\bn_{1\theta}\right),
\\
&&\bB_1 = \bnz\times\left(\bk\times\bn_{1\theta}\right),
\\
&&A_2 = \bnz\cdot\left(\bk\times\bn_{2\theta}\right) + \bn_1\cdot\left(\bk\times \bn_{1\theta}\right)
+ \bnz\cdot\curl\bnz,
\\
&&\bB_2 = \bnz\times\left(\bk\times\bn_{2\theta}\right) + \bn_1\times\left(\bk\times \bn_{1\theta}\right)
+ \bnz\times\curl\bnz,
\\
&&A_3 = \bnz\cdot\left(\bk\times\bn_{3\theta}\right) + \bn_1\cdot\left(\bk\times \bn_{2\theta}\right)
+  \bn_2\cdot\left(\bk\times \bn_{1\theta}\right)
\\
&&\qquad\qquad\qquad\qquad\qquad
+ \bnz\cdot\curl\bn_{1} + \bn_1\cdot\curl \bnz,
\\
&&\bB_3 = \bnz\times\left(\bk\times\bn_{3\theta}\right) + \bn_1\times\left(\bk\times \bn_{2\theta}\right)
+  \bn_2\times\left(\bk\times \bn_{1\theta}\right)
\\
&&\qquad\qquad\qquad\qquad\qquad
+ \bnz\times\curl\bn_{1} + \bn_1\times\curl \bnz.
\eeaa
We write these expressions as in \eq{Ahat}.

Using these expansions in \eq{dirPDE}, Taylor expanding the result,
and equating coefficients of powers of $\eps^{1/2}$, we get at the order $\eps^{-3/2}$
the linearized equations
\eq{linsysasy}. From Proposition~\ref{prop:linsys}, this system
has the twist-wave eigenvalue
$\omega^2=\Wtwist^2(\bk;\bnz)$, where $\Wtwist$ is given by \eq{bmode}.
The corresponding solution for $\{\bn_1,\lambda_1\}$ is
\be
\bn_1(\theta,\bx,t) = v(\theta,\bx,t) \Rtwist(\bx,t),\qquad
\lambda_1=0,
\label{soln1}
\ee
where $v$ is an arbitrary scalar-valued function,
and $\Rtwist$ is defined in \eq{Rtwist}.

Equating coefficients of the order $\eps^{-1}$, we obtain that
\be
\mathcal{L} \bn_{2\theta\theta} - \lambda_2\bnz = \bF_1,
\qquad
\bnz\cdot \bn_{2} = G_1,
\label{twistlin1}
\ee
where $\mathcal{L}$ is defined in \eq{defL}
and
\beaa
&&\bF_1 = -\beta\left[\rA_{2\theta}\left(\bk\times\bnz\right) + \bk\times\left(A_1\bn_1\right)_\theta
+ A_1\left(\bk\times\bn_{1\theta}\right)\right]
\\
&&\qquad
- \gamma \left[\bk\times\left(\rB_{2\theta}\times\bnz\right) + \bk\times\left(\bB_1\times\bn_1\right)_\theta
- \bB_1\times\left(\bk\times\bn_{1\theta}\right)\right] + \lambda_1 \bn_1,
\\
&&G_1 = -\frac{1}{2} \bn_1\cdot\bn_1.
\eeaa
Using Proposition~\ref{prop:linsys} to solve \eq{twistlin1}, we find that
\bea
&&\bn_2 =- r u \left(\bk\cdot\bnz\right) \Rsplay +  \frac{1}{2}v^2 (\bk\times\Rtwist) ,
\label{soln2}
\\
&&\lambda_2 = v_\theta^2 \left\{
\beta k^2 + r\left(\beta-\gamma\right)\left(\bk\cdot\bnz\right)^2\right\}
\left[k^2-\left(\bk\cdot\bnz\right)^2\right],
\nonumber
\eea
where $u(\theta,\bx,t)$ satisfies \eq{veq}, and
\[
r =  \frac{\beta - \gamma}{\alpha-\beta}.
\]
We could add an arbitrary scalar multiple of the null-vector $\Rtwist$ to
$\bn_2$, but this would not alter our final equations, so we omit it for simplicity.

Equating coefficients of the order $\eps^{-1/2}$, we obtain that
\be
\mathcal{L} \bn_{3\theta\theta} - \lambda_3\bnz = \bF_2,
\qquad
\bnz\cdot \bn_{3} = G_2,
\label{twistlin2}
\ee
where
\beaa
&&\bF_2 = 2\omega \bn_{1\theta t} + \omega_t \bn_{1\theta} + \lambda_2\bn_1 + \lambda_1 \bn_2
+ \alpha\left[\div \left(\bn_{1\theta} \right)\bk + \nabla\left(\bk\cdot\bn_{1\theta}\right)\right]
\\
&&\qquad\quad-\beta\left[ \rA_{3\theta}\left(\bk\times\bnz\right) + \bk\times\left(A_2\bn_1\right)_\theta
+ \bk\times\left(A_1\bn_2\right)_\theta + A_2\left(\bk\times \bn_{1\theta}\right)\right.
\\
&&\qquad\qquad\qquad\quad\left.
+ A_1 \left(\bk\times \bn_{2\theta}\right) + A_1\curl \bnz + \curl\left(A_1\bnz\right)\right]
\\
&&\qquad\quad- \gamma\left[ \bk\times\left(\rB_{3\theta}\times\bnz\right) + \bk\times\left(\bB_2\times\bn_1\right)_\theta
+ \bk\times\left(\bB_1\times\bn_2\right)_\theta \right.
\\
&&\qquad\qquad\qquad\quad- \bB_2\times\left(\bk\times \bn_{1\theta}\right)
- \bB_1\times \left(\bk\times \bn_{2\theta}\right)
\\
&&\qquad\qquad\qquad\qquad\quad\left.- \bB_1\times\curl \bnz + \curl\left(\bB_1\times\bnz\right)\right],
\\
&&G_2 = -\bn_1\cdot\bn_2.
\eeaa
From Proposition~\ref{prop:linsys}, this system is solvable if
$\Rtwist\cdot\bF_2 = 0$.
After some algebra, we compute that this solvability condition gives
\eq{ueq}.


\subsection{Weakly nonlinear polarized waves}
\label{app:pol}

We assume that $\bnz$ and $\bk$ are constant and $\bk = k\bnz$
is parallel to $\bnz$. We look for an expansion of the same form
as the one used in the previous section for twist waves,
with a plane-wave phase $\Phi(x,t) = \bk\cdot \bx - \omega t$.

When $\bk$ is parallel to $\bnz$,
the linearized dispersion relation \eq{lindisp} reduces to
$\left[\omega^2 - \gamma k^2\right]^2 = 0$, so $\omega^2=\gamma k^2$.
From Proposition~\ref{prop:linsys},
the solution of the leading-order $O(\eps^{-3/2})$-equation \eq{linsysasy}
is then
\[
\bn_1(\theta,\bx,t) = \bu(\theta,\bx,t),\qquad \lambda_1=0,
\]
where $\bu$ is an arbitrary vector such that $\bnz\cdot\bu = 0$.

The $O(\eps^{-1})$-equation \eq{twistlin1} simplifies to
\beaa
&&\mathcal{L} \bn_{2\theta\theta} - \lambda_2 \bnz
= -\gamma k^2 \left(\bu_\theta\cdot\bu_\theta\right) \bnz,
\\
&&\bnz\cdot \bn_2 = -\frac{1}{2} (\bu\cdot\bu),
\eeaa
whose solution is
\beaa
&&\bn_2 = -\frac{1}{2} (\bu\cdot\bu) \bnz,
\\
&&\lambda_2 = \frac{1}{2}(\alpha-\gamma)  k^2 \left(\bu\cdot\bu\right)_{\theta\theta}
+ \gamma k^2 \left(\bu_\theta\cdot\bu_\theta\right).
\eeaa
We could add to $\bn_2$ an arbitrary vector orthogonal to $\bnz$,
but this would not alter our final equations, so we omit it for simplicity.

Computing $\bF_2$ in the $O(\eps^{-1/2})$-equation \eq{twistlin2}
and imposing the solvability condition in Proposition~\ref{prop:linsys} that $\bF_2$ is parallel
to $\bnz$, we get after some algebra that
\bea
&&\bu_{\theta t} + \frac{\omega}{k}\bnz\cdot\nabla \bu_\theta
+ \frac{\alpha k^2}{2\omega}\left( \bu\cdot\bu_\theta\right)_{\theta} \bu
\nonumber
\\
&&\qquad- \frac{\beta k^2}{2\omega}\left\{
\left[\bu\cdot\left(\bnz\times\bu_{\theta\theta}\right)\right](\bnz\times\bu)
+ 2\left[\bu\cdot\left(\bnz\times\bu_\theta\right)\right]
\left(\bnz\times\bu_\theta\right)\right\}
\nonumber
\\
&&\qquad- \frac{\gamma k^2}{2\omega}
\left\{\left[\left(\bu\cdot\bu\right) \bu_{\theta}\right]_\theta
- \left(\bu_\theta\cdot\bu_\theta\right) \bu\right\}
= 0.
\label{poleq1}
\eea
A computation in terms of components show that if $\bu$ is orthogonal to
the constant unit vector $\bnz$, then
\beaa
&&\left[\bu\cdot\left(\bnz\times\bu_{\theta\theta}\right)\right](\bnz\times\bu)
+ 2\left[\bu\cdot\left(\bnz\times\bu_\theta\right)\right]\left(\bnz\times\bu_\theta\right)
\\
&&\qquad\qquad\qquad
= \left(\bu\cdot\bu_\theta\right)_\theta \bu - \left[\left(\bu\cdot\bu\right)\bu_\theta\right]_\theta
+\left(\bu_\theta\cdot\bu_\theta\right)\bu.
\eeaa
Using this result in \eq{poleq1}, we obtain \eq{poleq}.



\begin{thebibliography}{99}

\bibitem{saxton} R.~Saxton, Dynamic instability of the liquid crystal director,
in \emph{Contemp. Math.}, Vol. 100, Current Progress in Hyperbolic Systems,
W. B. Lindquist ed., AMS, Providence, RI, 1989, 325--330.


\bibitem{hs} J.~K.~Hunter, and R.~Saxton,
Dynamics of director fields, {\em SIAM J.~Appl.~Math.},
\textbf{51} (1991), 1498-1521.


\bibitem{ghz} R. Glassey, J. K. Hunter, and Y. Zheng, Singularities
of a variational wave equation, \emph{J. Diff. Eq.}, \textbf{129} (1996), 49--78.


\bibitem{dafermos} C.~M.~Dafermos, \emph{Hyperbolic Conservation Laws in Continuum Physics},
2nd ed., Springer-Verlag, Berlin, 2005.


\bibitem{bz} A.~Bressan, and Y.~Zheng, Conservative solutions to a nonlinear variational
wave equation, to appear in \emph{Comm. Math. Phys.}.

\bibitem{hz}  J.~K.~Hunter, and  Y.~Zheng, On a completely integrable
nonlinear hyperbolic variational equation, \emph{Physica D}, \textbf{79},
(1994), 361--386.

\bibitem{bs} R. Beals, D. Sattinger, and J. Szmigielske,
Inverse scattering solutions of the Hunter-Saxton equation, \emph{Applicable Analysis},
\textbf{78} (2001), 255--269.

\bibitem{km} B.~Khesin and G. Misiolek, Euler equations on homogeneous spaces
and Virasoro orbits, \emph{Adv. Math}, \textbf {176} (2003), 116--144.

\bibitem{hz1}  J.~K.~Hunter, and  Y.~Zheng,
On a nonlinear hyperbolic variational wave equation I. Global existence of
weak solutions, \emph{Arch. Rat. Mech. Anal.}, \textbf{129} (1995), 355--383.

\bibitem{zz} P.~Zhang, and Y.~Zheng, Existence and uniqueness of solutions of an asymptotic
equation arising from a variational wave equation with general data,
\emph{Arch. Rat. Mech. Anal.}, \textbf{155} (2000), 49--83.

\bibitem{bc} A.~Bressan, and A.~Constantin, Global solutions of the Hunter-Saxton
equation, \emph{SIAM J. Math. Anal.},
\textbf{37} (2005), 996-1026.

\bibitem{CH} R.~Camassa and D.~Holm, An integrable shallow water equation with peaked solitons,
\emph{Phys. Rev. Lett.}, \textbf{71} (1993), R51--R79.

\bibitem{zz1} A. Bressan, P.~Zhang, and Y.~Zheng, On asymptotic variational wave equations,
to appear in \emph{Arch. Rat. Mech. Anal.}.

\bibitem{deG} P.~G.~de Gennes, and J.~Prost, \textit{The Physics of Liquid Crystals},
2nd Edition, Oxford University Press, Oxford, 1993.

\bibitem{ah} G.~Al\`\i, and J.~K.~Hunter, Linearly degenerate variational
wave equations, in preparation.

\bibitem{cmp} G. F. Mazenko, \emph{Fluctuations, Order, and Defects},
John Wiley \& Sons, New Jersey, 2003.

\bibitem{ah1} G.~Al\`\i, G., and J.~K.~Hunter, Diffractive nonlinear geometrical optics
for variational wave equations and the Einstein equations,
to appear in \emph{Comm. Pure Appl. Math,}

\bibitem{ericksen} J. L. Ericksen, Twisting of liquid crystals,
\emph{J. Fluid Mech.}, \textbf{27} (1967), 59--64.

\bibitem{ericksen1} J. L. Ericksen,
Twist waves in liquid crystals,
\emph{Q. J. Mech. Appl. Math.} (1972), \textbf{25}, 463--465.

\bibitem{shahinpoor} M. Shahinpoor,
Finite twist waves in liquid crystals,
\emph{Q. J. Mech. Appl. Math.}, \textbf{28} (1975), 223--231.

\bibitem{bh} M.~Brio and J. K. Hunter, Rotationally invariant hyperbolic waves,
\emph{Comm. Pure Appl. Math.}, \textbf{43} (1990), 1037--1053.

\bibitem{per} R. Glassey, J. K. Hunter, and Y. Zheng, Singularities
and oscillations in a nonlinear variational wave equation, in
\emph{Singularities and Oscillations}, ed. J. Rauch and M. Taylor,
IMA Volumes in Mathematics and its Applications, \textbf{91}, Springer-Verlag,
1996, 37--60.

\bibitem{rev} J.~K.~Hunter,
Asymptotic equations for nonlinear hyperbolic waves,
in {\em Surveys in Applied Mathematics}, Vol. 2,
ed. M. Freidlin et.al., 167--276, Plenum Press, 1995.
\end{thebibliography}
\end{document}